\documentclass{birkjour}

\usepackage{latexsym}
\usepackage{amssymb}
\usepackage{amsmath}
\usepackage{amscd}
\usepackage{mathrsfs}
\newtheorem{Th}{Theorem}[section]
\newtheorem{Co}[Th]{Corollary}

\newtheorem{Lem}[Th]{Lemma}

\newtheorem{Rem}[Th]{Remark}

\newtheorem{Pro}[Th]{Proposition}
\newcommand{\demo}{\par\noindent{\it Proof. \/}\ }
\newcommand{\enD}{\hfill $\Box$ \vspace{3truemm}\par}
\newcommand{\bx}{\mbox{\boldmath $x$}}
\newcommand{\bX}{\mbox{\boldmath $X$}}
\newcommand{\be}{\mbox{\boldmath $e$}}
\newcommand{\ba}{\mbox{\boldmath $a$}}
\newcommand{\bb}{\mbox{\boldmath $b$}}

\newcommand{\bv}{\mbox{\boldmath $v$}}

\newcommand{\bo}{\mbox{\boldmath $0$}}

\newcommand{\bw}{\mbox{\boldmath $w$}}

\newcommand{\bn}{\mbox{\boldmath $n$}}

\newcommand{\R}{{\mathbb R}}
\newcommand{\lon}{\longrightarrow}
\newcommand{\vect}[1]{\boldsymbol{#1}}
\newcommand{\ve}{\varepsilon}
\newcommand{\inner}[2]{\left\langle{#1},{#2}\right\rangle_1}

\newcommand{\de}{\partial_{\varepsilon}}
\newcommand{\deo}{ \left. \partial_{\varepsilon} \right|_{\varepsilon=0}}
\newcommand{\Pde}{\frac{\partial}{\partial \varepsilon}}
\newcommand{\Pdeo}{ \left. \frac{\partial}{\partial \varepsilon} \right|_{0}}
\newcommand{\bl}{\mbox{\boldmath $\ell$}}

\newcommand{\hpab}{h_{ab}(\bn^{T},\pm\bn^{S})}
\newcommand{\Wpm}[2]{h^{#1}_{#2}(\bn^{T},\pm\bn^{S})}

\newcommand{\hpmij}{h_{ij}(\bn^{T},\pm\bn^{S})}
\newcommand{\hpij}{h_{ij}(\bn^{T},\bn^{S})}

\newcommand{\Hpm}{H_{\ell}(\bn^{T},\pm\bn^{S})}
\newcommand{\Hp}{H_{\ell}(\bn^{T},\bn^{S})}

\newcommand{\Kpm}{K_{\ell}(\bn^{T},\pm\bn^{S})}
\newcommand{\Kp}{K_{\ell}(\bn^{T},\bn^{S})}

\begin{document}
\title[Marginally trapped surfaces in Minkowski space-time]{The lightlike geometry of  marginally \\ trapped surfaces in
Minkowski space-time}
\author[A. Honda]{Atsufumi HONDA}
\address{%
National Insitute of Technology, 
Miyakonojo College \\
Miyakonojo 885-8567\\ 
Japan}

\email{atsufumi@cc.miyakonojo-nct.ac.jp
}

\thanks{%
Partly supported by the Grant-in-Aid for JSPS Fellows, 
The Ministry of Education, Culture, Sports, Science and Technology, Japan.
}
\author[S. Izumiya]{Shyuichi IZUMIYA} 
\address{%
Department of Mathematics
\\ Hokkaido University
\\ Sapporo 060-0810
\\ Japan} 

\email{izumiya@math.sci.hokudai.ac.jp}

%\thanks{%
%Work partially supported  by Grant-in-Aid for Scientific Research (No. 21654007, No. 22340011) 
%Japan Society for the Promotion of Science.}
\date{\today}
\begin{abstract}
The lightlike geometry of codimension two spacelike submanifolds in Lorentz-Minkowski space has been developed in
\cite{IR07} which is a natural Lorentzian analogue of the classical Euclidean differential geometry of hypersurfaces.
In this paper we investigate a special class of surfaces (i.e., {\it marginally trapped surfaces}) in Minkowski space-time
from the view point of the lightlike geometry.

\end{abstract}

\maketitle

% \tableofcontents
\section{Introduction}

The notion of trapped surfaces introduced by Roger Penrose \cite{Pen65}
is one of the most important subjects in cosmology and general relativity.
It plays a principal role for the singularity theorems, the analysis of gravitational
collapse, the cosmic censorship hypothesis, Penrose inequality, etc.
A marginally trapped surface is
defined to be a spacelike surface with isotropic (i.e., zero or lightlike) mean curvature vector field, which separates the trapped surfaces from the untrapped one.
The surface of a black hole might be the marginally trapped surface.
Mathematically, marginally trapped surfaces are viewed as space-time analogues of minimal surfaces in Riemannian manifolds.
In Minkowski space-time (i.e., Lorentz-Minkowski $4$-space),
the spacelike surface with the zero mean curvature vector field is said to be minimal. However, we call it strongly marginally trapped in this paper.
In \cite{AGM05} the Weierstrass-Bryant type representation of a marginally trapped surface in Minkowski space-time was given.
\par
On the other hand, it has been claimed in \cite{IR07} that codimension two spacelike submanifolds in Lorentz-Minkowski space play a role similar to hypersurfaces in Euclidean space.
For a codimension two spacelike submanifold in Lorentz-Minkowski space, we have two lightlike normal directions along the submanifolds,
if we use these lightlike normal directions as the unit normal vector field along a hypersurface in Euclidean space, we can
define curvatures, etc. This is the basic idea of the lightlike geometry for codimension two spacelike submanifolds\cite{IKPR06,IR07,IRS09}.
In this paper, inspired by the beautiful survey article\cite{CGP10}, we apply the lightlike geometry to marginally trapped surfaces in Minkowski space-time.
Then we observe natural geometric properties analogous to minimal surfaces in Euclidean $3$-space.
For example, a marginally trapped surface is a spacelike surface such that one of the lightcone mean curvatures vanishes (Proposition 3.2 and Corollary 4.4).
One of the consequence of the lightlike geometry is that a totally umbilical marginally trapped surface is given by a graph of 
a smooth function $f(u_1,u_2)$ (Corollary 3.5).
Therefore, there are so many complete marginally trapped surfaces in Minkowski space-time compared with
minimal surfaces in Euclidean space. 
Motivated by this fact, we obtain a partial differential equation for a marginally trapped graph surface of
two functions $f(u_1,u_2),g(u_1,u_2)$ (Theorem 6.1).
As a consequence, a graph surface of a single function $f(u_1,u_2)$ is strongly marginally trapped if and only if
$f$ is harmonic.
Therefore, the Bernstein theorem\cite{Be} does not hold even for strongly marginally trapped surfaces.
We also have a characterization of marginally trapped surfaces by the variational problem of the area functionals with respect to
lightlike normal directions (Theorem 5.2).
We remark that the class of marginally trapped surfaces includes a generalization of the notion of
not only  minimal surfaces in Euclidean $3$-space, maximal surfaces in Lorentz-Minkowski $3$-space,
the surfaces with constant mean curvature $\pm 1$ (briefly, CMC$\pm 1$ surfaces) in Hyperbolic $3$-space and CMC$\pm 1$ spacelike surfaces in de Sitter $3$-space but also intrinsic flat spacelike surfaces in the 
lightcone (Theorem 7.3).
\par
In \S 2 we explain the basic facts and notations of Minkowski space-time. 
Basic properties of the lightlike geometry are given in \S 3 and \S 4.
We consider the variation problem with respect to lightlike normal directions in \S 5.
In \S 6, we consider partial differential equations which are satisfied by marginally trapped or strongly marginally trapped surfaces. In \S 7 we explain how naturally interpret minimal surfaces in Euclidean $3$-space, maximal surfaces in Lorentz-Minkowski $3$-space,
CMC$\pm 1$ surfaces in Hyperbolic $3$-space, CMC$\pm 1$ spacelike surfaces in de Sitter $3$-space and intrinsic flat spacelike surfaces in the 
lightcone as marginally trapped surfaces.

\section{Basic facts and notations in Minkowski space-time}
\par
In this section we briefly introduce basic notations on Minkowski space-time.
The detailed properties should be referred in \cite{Oneil}.
Denote by $\mathbb{R}^{4}_1$ the Lorentz-Minkowski $4$-dimensional space $(\mathbb{R}^
{4}, \langle , \rangle_1)$, with the pseudo scalar product given by
$$\langle \mbox{\boldmath$x$},\mbox{\boldmath$y$}\rangle _1 =
-x_0y_0+x_1y_1+x_2y_2+x_3y_3,$$
where $\mbox{\boldmath$x$}=(x_0,x_1,x_2,x_3),\mbox{\boldmath$x$}=(y_0,y_1,y_2,y_3).$ 
The {\it norm} of a vector $\mbox{\boldmath$x$}$ is defined to be
$\|\mbox{\boldmath$x$}\|_1=\sqrt{|\langle \mbox{\boldmath$x$},\mbox{\boldmath$x$}\rangle_1|}$.
A vector $\mbox{\boldmath$x$} = (x_0,x_1,x_2,x_3) \in \mathbb{R}^{4}_1 \backslash \{
\mbox{\boldmath$0$} \}$ is said to be {\it  spacelike, timelike} or {\it
lightlike} according to $\langle \mbox{\boldmath$x$}, \mbox{\boldmath$x$}\rangle_1 >0$, $=0$,
or $<0$, respectively. For any $\mbox{\boldmath$v$} \in \mathbb{R}_1^{4}\setminus  \{\mbox{\boldmath$0$}\}$ and $c\in
\mathbb{R}$, the hyperplane  with the pseudo-normal \mbox{\boldmath$v$} is given by
$$HP(\mbox{\boldmath$v$},c)=\{\mbox{\boldmath$v$} \in \mathbb{R}_1^{4} \vert \langle \mbox{\boldmath$x$}, \mbox{\boldmath$v$}
\rangle _1=c\}.$$ We say that  $HP(\mbox{\boldmath$v$},c)$ is {\it spacelike,
timelike} or {\it lightlike} provided $\mbox{\boldmath$v$}$ is timelike,
spacelike or lightlike respectively.

The {\it Hyperbolic $n$-space} is defined in
this context as the subset
$$ \;\; H^3_+(-1)=\{\mbox{\boldmath$x$}\in \mathbb{R}^{4}_1 |
\langle \mbox{\boldmath$x$} ,\mbox{\boldmath$x$}\rangle_1 =-1, \; x_0>0 \}. $$
Other well known
hypersurfaces in Lorentz-Minkowski space are the {\it de Sitter $3$-space:}
$$
S^3_1=\{\mbox{\boldmath$x$}\in \mathbb{R}^{4}_1 | \langle \mbox{\boldmath$x$} ,\mbox{\boldmath$x$}\rangle_1 =1\ \}, $$
and
the {\it open lightcone (at the origin):}
$$
LC^*=\{\mbox{\boldmath$x$}\in \mathbb{R}^{4}_1  \ | \mbox{\boldmath$x$}\not= \mbox{\boldmath$0$},\
\langle \mbox{\boldmath$x$} ,\mbox{\boldmath$x$} \rangle_1 =0\}. $$
We shall also consider the {\it future lightcone:}
$$LC^*_+=\{ \mbox{\boldmath$x$} \in LC^*\ | x_0> 0\}
$$
and the {\it lightcone unit $2$-sphere:} $$
S^{2}_+=\{\mbox{\boldmath$x$} =(x_0,x_1,x_2,x_3)\ |\ \langle \mbox{\boldmath$x$} ,\mbox{\boldmath$x$} \rangle_1
=0,\, x_0=1\}.$$
 If  $\mbox{\boldmath$x$} =(x_0,x_1,x_2
,x_3)$ is a non-zero lightlike vector, we have $x_0 \not= 0$ and
write
$$\widetilde{\mbox{\boldmath$x$}}=\left(1,\frac{x_1}{x_0},\frac{x_2}{x_0}
,\frac{x_3}{x_0}\right)\in S^{2}_+.$$
Therefore, we define a projection $\pi ^L_S:LC^*\lon S^2_+$ by $\pi^L_S(\mbox{\boldmath$x$})=\widetilde{\mbox{\boldmath$x$}}.$

If  $\{\mbox{\boldmath$e$}_{0},\mbox{\boldmath$e$}_{1},\mbox{\boldmath$e$}_{2}  ,
\mbox{\boldmath$e$}_{3}\}$ represents the canonical
basis of $\mathbb{R}^{4}_1$, we shall say that a timelike vector $\mbox{\boldmath$v$}$ is
{\it future directed} if $\langle \mbox{\boldmath$v$} ,\mbox{\boldmath$e$} _0\rangle_1 <0.$

Given $3$ vectors $\mbox{\boldmath$x$}_1,\mbox{\boldmath$x$}_2,\mbox{\boldmath$x$}_3  \in \mathbb{R}^{4}_1,$
we can consider their  wedge product,
\[
\mbox{\boldmath$x$}_1 \wedge \mbox{\boldmath$x$}_2 \wedge \mbox{\boldmath$x$}_3= \left|
\begin{array}{cccc}
-\mbox{\boldmath$e$}_{0}&\mbox{\boldmath$e$}_{1}& \mbox{\boldmath$e$}_2 &\mbox{\boldmath$e$}_{3}\\
x^1_{0}& x^1_{1}& x^1_{2} &x^1_{3}\\
x^2_{0}& x^2_{1}&  x^2_{2} & x^2_{3}\\
x^3_{0}&x^3_{1}& x^3_{2} & x^3_{3}
\end{array}
\right| ,
\]
where  $\mbox{\boldmath$x$}_i=(x_0^i,x_1^i, x_2^i ,x_3^i), \; i=1,2,3.$

Clearly $ \langle \mbox{\boldmath$x$},\mbox{\boldmath$x$}_1\wedge \mbox{\boldmath$x$}_2 \wedge 
\mbox{\boldmath$x$}_3  \rangle_1 ={\rm
det}(\mbox{\boldmath$x$},\mbox{\boldmath$x$}_1,\mbox{\boldmath$x$}_2, \mbox{\boldmath$x$}_3). $ So $\mbox{\boldmath$x$}_1\wedge\mbox{\boldmath$x$}_2\wedge\mbox{\boldmath$x$}_3$ is pseudo orthogonal to $\mbox{\boldmath$x$}_i, i=1,2,3.$
\section{Lightcone curvatures of spacelike surfaces in Minkowski space-time}

 We introduce the basic geometrical tools for the study of spacelike surfaces in
Minkowski space-time. 
Consider the orientation of $\mathbb{R}^{4}_1$
provided by the volume form  $\mbox{\boldmath$l$}_0 \wedge \cdots \wedge \mbox{\boldmath$l$}_4$,
where $\{\mbox{\boldmath$l$}_i\}_{i=0}^4$ is the dual basis of the canonical basis
$\{\mbox{\boldmath$e$}_i\}_{i=0}^4$. We also give  $\mathbb{R}^{4}_1$ a timelike
orientation by choosing $\mbox{\boldmath$e$} _0=(1,0,0,0)$ as a future timelike
vector field.

Given a surface $M$  in  $\mathbb{R}^{4}_1$ consider a local  parametrization (embedding) $\mbox{\boldmath$X$}:U\rightarrow
\mathbb{R}^{4}_1$ of $M$, where  $U$ is an open subset of $\mathbb{R}^{2}.$ We write
$M=\mbox{\boldmath$X$}(U)$ and identify $M$ and $U$ through the embedding $\mbox{\boldmath$X$}.$ We
say that $\mbox{\boldmath$X$}$ is {\it spacelike} if $\mbox{\boldmath$X$}_{u_i}$ $i=1,2$
are always spacelike vectors. Therefore, the tangent space $T_p  M$
of $M$ at $p$ is a spacelike subspace (i.e., consists of spacelike
vectors) for any point $p\in M$. In this case, the pseudo-normal
space $N_p M$ is a timelike plane (i.e., Lorentz plane).
We denote by $N(M)$ the pseudo-normal bundle
over $M.$ We can arbitrarily choose a
future directed unit timelike normal section $\mbox{\boldmath$n$} ^T(u)\in N_p(M)$,
where $p=\mbox{\boldmath$X$}(u).$ 
We remark that $\bn^T$ always exists even globally.
Here, we say that $\mbox{\boldmath$n$}^T$ is {\it future directed}
if $\langle \mbox{\boldmath$n$}^T ,{\mbox{\boldmath$e$}} _0\rangle_1 <0.$
Therefore we can construct  a
spacelike unit normal section $\mbox{\boldmath$n$}^S(u)\in N_p(M)$ by
\begin{equation}
\mbox{\boldmath$n$} ^S(u)=\frac{\mbox{\boldmath$n$}^T (u )\wedge\mbox{\boldmath$X$} _{u_1}(u )\wedge
\mbox{\boldmath$X$} _{u_{2}}(u )}{\|\mbox{\boldmath$n$}^T (u )\wedge \mbox{\boldmath$X$}_{u_1}(u )
\wedge \mbox{\boldmath$X$}_{u_{2}}(u )\|_1},
\end{equation}
and we have
$$\langle \mbox{\boldmath$n$}^T, \mbox{\boldmath$n$}^T\rangle_1 =-1,  \;\;\; \langle
\mbox{\boldmath$n$}^T,\mbox{\boldmath$n$}^S\rangle_1 =0, \;\;\; \langle \mbox{\boldmath$n$}^S,\mbox{\boldmath$n$}^S\rangle_1 =1.$$
Although we
could also choose $-\mbox{\boldmath$n$} ^S(u)$ as a spacelike unit normal section
with the above properties, we shall fix the direction $\mbox{\boldmath$n$} ^S(u)$
throughout this section.
In the global sense, $\bn^S$ exists for an orientable surface $M.$
We call $(\mbox{\boldmath$n$} ^T,\mbox{\boldmath$n$}^S)$ a {\it future
directed normal frame\/} along $M=\mbox{\boldmath$X$}(U).$ Clearly, the vectors
$\mbox{\boldmath$n$}^T (u)\pm \mbox{\boldmath$n$}^S(u)$ are lightlike. We choose here $\mbox{\boldmath$n$} ^T+\mbox{\boldmath$n$}^S$
as a lightlike normal vector field along $M.$ Since
$\{\mbox{\boldmath$X$}_{u_1}(u),\mbox{\boldmath$X$}_{u_{2}}(u)\}$ is a basis of $T_pM,$ the
vectors $ \{\mbox{\boldmath$n$}^T (u),\mbox{\boldmath$n$} ^S(u),\mbox{\boldmath$X$}_{u_1}(u)
,\mbox{\boldmath$X$}_{u_{2}}(u)\} $ provide a basis for $T_p \mathbb{R}^{4}_1.$
In \cite{IR07} it has been shown the following lemma.

\begin{Lem}
\label{l31}  For given two future directed unit timelike normal sections \newline
$\mbox{\boldmath$n$}^T(u), \overline{\mbox{\boldmath$n$}}^T(u)\in N_p(M),$  the corresponding lightlike
normal sections \newline $\mbox{\boldmath$n$}^T(u)\pm\mbox{\boldmath$n$}^S(u), \overline{\mbox{\boldmath$n$}}^T(u)\pm\overline{\mbox{\boldmath$n$}}^S(u)$
are parallel.It follows that 
\[
\pi ^L_S(\mbox{\boldmath$n$}^T(u)\pm\mbox{\boldmath$n$}^S(u))=\pi^L_S(\overline{\mbox{\boldmath$n$}}^T(u)\pm\overline{\mbox{\boldmath$n$}}^S(u))
\]
\end{Lem}

Under the identification of $M$ and $U$ through $\mbox{\boldmath$X$},$ we have the
linear mapping given by the derivative of the lightcone normal
vector field $\mbox{\boldmath$n$}^T\pm\mbox{\boldmath$n$}^S$ at each point $p=\mbox{\boldmath$X$}(u)\in M$,
$$
d_p(\mbox{\boldmath$n$}^T\pm\mbox{\boldmath$n$}^S):T_pM\rightarrow T_p \mathbb{R}^{4}_1= T_pM\oplus N_p(M).$$
Consider the orthogonal projections 
$$
\pi ^t:T_pM \oplus
N_p(M)\rightarrow T_p(M)$$ and $$\pi ^n:T_p(M)\oplus N_p(M)\rightarrow
N_p(M).
$$ 
We define the {\it $(\mbox{\boldmath$n$}^T,\pm\mbox{\boldmath$n$}^S)$-shape operator} of $M=\mbox{\boldmath$X$} (U)$
at $p=\mbox{\boldmath$X$} (u)$
\begin{equation}
-d_p(\mbox{\boldmath$n$} ^T\pm\mbox{\boldmath$n$}^S)^t=-\pi ^t\circ d_p(\mbox{\boldmath$n$}^T\pm\mbox{\boldmath$n$}^S)
\end{equation}
and denote it by $S_{p}(\mbox{\boldmath$n$}^T,\pm\mbox{\boldmath$n$}^S)$.

The {\it normal connection with respect to
$(\mbox{\boldmath$n$}^T,\pm\mbox{\boldmath$n$}^S)$} of $M$ at $p$ is defined by the linear transformation
\begin{equation}
d_p(\mbox{\boldmath$n$} ^T\pm\mbox{\boldmath$n$}^S)^n=\pi ^n\circ d_p(\mbox{\boldmath$n$}^T\pm\mbox{\boldmath$n$}^S).
\end{equation}
We also define
$S_p(\mbox{\boldmath$n$}^T)=-\pi^t\circ d_p\mbox{\boldmath$n$}^T$ and
$S_p(\mbox{\boldmath$n$}^S)=-\pi\circ d_p\mbox{\boldmath$n$}^S.$
We respectively call these a {\it $\mbox{\boldmath$n$}^T$-shape operator} and a {\it
$\mbox{\boldmath$n$}^S$-shape operator} of $M=\mbox{\boldmath$X$}(U)$ at $p=\mbox{\boldmath$X$}(u).$
The eigenvalues of $S_{p}(\mbox{\boldmath$n$}^T,\pm\mbox{\boldmath$n$}^S)$, which is denoted by $\{\kappa
_{i}(\mbox{\boldmath$n$}^T,\pm\mbox{\boldmath$n$}^S)(p)\}_ {i=1}^{2},$ are called the {\it lightcone
principal curvatures \/}  with respect to $(\mbox{\boldmath$n$}^T,\pm\mbox{\boldmath$n$}^S) $ at $p$.
We also define the {\it $\mbox{\boldmath$n$}^T$-principal curvature} $\kappa _i(\mbox{\boldmath$n$}^T)(p)$
(respectively, {\it $\mbox{\boldmath$n$}^S$-principal curvature} $\kappa _i(\mbox{\boldmath$n$}^S)(p)$)
as the eigenvalues of $S_p(\mbox{\boldmath$n$}^T)$ (respectively, {\it $\mbox{\boldmath$n$}^T$-principal curvature} $\kappa _i(\mbox{\boldmath$n$}^S)(p)$).
Since 
$$S_{p}(\mbox{\boldmath$n$}^T,\pm\mbox{\boldmath$n$}^S)=S_p(\mbox{\boldmath$n$}^T)\pm S_p(\mbox{\boldmath$n$}^S),
$$ we have $\kappa
_{i}(\mbox{\boldmath$n$}^T,\pm\mbox{\boldmath$n$}^S)(p)=\kappa _i(\mbox{\boldmath$n$}^T)(p)\pm \kappa_i(\mbox{\boldmath$n$}^S)(p).$
%We define the notion of curvature as follows.
The {\it lightcone Gauss-Kronecker curvature} with respect to
$(\mbox{\boldmath$n$}^T,\pm\mbox{\boldmath$n$}^S)$ at $p$ is defined as follows:
\begin{equation}
K_\ell(\mbox{\boldmath$n$}^T,\pm\mbox{\boldmath$n$}^S)(p)={\rm det} S_{p}(\mbox{\boldmath$n$}^T,\pm\mbox{\boldmath$n$}^S).
\end{equation}
The {\it lightcone mean curvature} with respect to $(\mbox{\boldmath$n$}^T,\pm\mbox{\boldmath$n$}^S)$ at $p$ is defined by
\[
H_\ell(\mbox{\boldmath$n$}^T,\pm\mbox{\boldmath$n$}^S)(p)
=\frac{1}{2}{\rm Trace}\, S_{p}(\mbox{\boldmath$n$}^T,\pm\mbox{\boldmath$n$}^S).
\]
\par
On the other hand, the {\it mean curvature vector} of $M=\mbox{\boldmath$X$}(U)$ at
$p=\mbox{\boldmath$X$}(u)$ is defined to be
\[
\mathfrak{H}(p)=\frac{1}{2}{\rm Trace}\,S_p(\mbox{\boldmath$n$}^T)\mbox{\boldmath$n$}^T(u)+\frac{1}{2}{\rm Trace}\,S_p(\mbox{\boldmath$n$}^S)\mbox{\boldmath$n$}^S(u).
\]
We have the following proposition.
\begin{Pro}
Under the above notations, the following conditions are equivalent.
\par
{\rm (1)} $H_\ell(\mbox{\boldmath$n$}^T,\mbox{\boldmath$n$}^S)(p)=0$ or $H_\ell(\mbox{\boldmath$n$}^T,-\mbox{\boldmath$n$}^S)(p)=0$.
\par
{\rm (2)} The mean curvature vector $\mathfrak{H}(p)$ is an isotropic vector, 
\par
where a vector is isotropic if it is a zero vector or a lightlike vector.
\end{Pro}

\demo
$S_{p}(\mbox{\boldmath$n$}^T,\pm\mbox{\boldmath$n$}^S)=S_p(\mbox{\boldmath$n$}^T)\pm S_p(\mbox{\boldmath$n$}^S),$ we have 
\[
{\rm Trace}\, S_{p}(\mbox{\boldmath$n$}^T,\pm\mbox{\boldmath$n$}^S)={\rm Trace}\, S_p(\mbox{\boldmath$n$}^T)\pm {\rm Trace}\, S_p(\mbox{\boldmath$n$}^S).
\]
On the other hand,
\[
4\langle \mathfrak{H}(p),\mathfrak{H}(p)\rangle _1=
-{\rm Trace}\, S_p(\mbox{\boldmath$n$}^T)^2+{\rm Trace}\, S_p(\mbox{\boldmath$n$}^S)^2.
\]
Thus, $\langle \mathfrak{H}(p),\mathfrak{H}(p)\rangle _1=0$ if and only if
${\rm Trace}\, S_p(\mbox{\boldmath$n$}^T)=\pm {\rm Trace}\, S_p(\mbox{\boldmath$n$}^S).$
This condition is equivalent to the condition (1).
\enD
We denote 
\[
H(\bn^T)(p)=({\rm Trace}\, S_p(\bn^T))/2\ \mbox{and}\ H(\bn^S)(p)=({\rm Trace}\, S_p(\bn^S))/2,
\]
which we respectively call the {\it mean curvature} of $M=\bX(U)$ at $p=\bX(u)$ {\it with respect
to} $\bn^T$ and $\bn^S$.
\par
On the other hand,
there is a concept of {\it trapped surfaces} in a space-time introduced by
Penrose in \cite{Pen65} which plays an extremely important role in 
cosmology and general relativity. In terms of mean curvature vector, a spacelike surface
in a space-time is {\it marginally trapped} if its mean curvature vector is isotropic at
each point.
The above proposition asserts that $M\subset \mathbb{R}^4_1$ is marginally trapped if
and only if $H_\ell (\mbox{\boldmath$n$}^T,\pm\mbox{\boldmath$n$}^S)\equiv 0$ for
any future directed normal frame $(\mbox{\boldmath$n$}^T,\mbox{\boldmath$n$}^S)$.
We also say that  $M=\bX(U)$ is a {\it strongly marginally trapped surface} if
the mean curvature vector $\mathfrak{H}$ is zero at any point. In some articles, this notion is
called {\it minimal}.
However, this class of surfaces includes both the minimal surfaces in Euclidean space
and the maximal surfaces in Minkowski $3$-space (cf., \S 6).
Therefore we call these strongly marginally trapped surfaces.
By the theorem on space-time singularities in \cite{Pen65}, there are no compact marginally trapped surfaces
in Minkowski space-time.
By Proposition 3.2, we have the following corollary.
\begin{Co} With the same notations as those in Proposition 3.2, we have the following{\rm :}
\par\noindent
{\rm (1)} $M$ is marginally trapped if and only if 
\[
H_\ell(\mbox{\boldmath$n$}^T,\mbox{\boldmath$n$}^S)\equiv 0\ \mbox{or}\ H_\ell(\mbox{\boldmath$n$}^T,-\mbox{\boldmath$n$}^S)\equiv 0.
\]
\par\noindent
{\rm (2)} $M$ is strongly marginally trapped if and only if 
\[
H_\ell(\mbox{\boldmath$n$}^T,\mbox{\boldmath$n$}^S)\equiv 0\ \mbox{and}\ H_\ell(\mbox{\boldmath$n$}^T,-\mbox{\boldmath$n$}^S)\equiv 0.
\]
\end{Co}
\par
On the other hand,
we say that a point $p$ is a {\it $(\mbox{\boldmath$n$}^T,\pm\mbox{\boldmath$n$}^S)$-umbilical point} if
all the principal curvatures coincide at $p$ and thus
$S_{p}(\mbox{\boldmath$n$}^T,\pm\mbox{\boldmath$n$}^S)=\kappa (\mbox{\boldmath$n$}^T,\pm\mbox{\boldmath$n$}^S)(p) 1_{T_{p}M}$, for some
function $\kappa.$ We say that $M$ is {\it totally
$(\mbox{\boldmath$n$}^T,\pm\mbox{\boldmath$n$}^S)$-umbilical} if all points on $M$ are
$(\mbox{\boldmath$n$}^T,\pm\mbox{\boldmath$n$}^S)$-umbilical.
In \cite{IR07} we have shown that a totally $(\mbox{\boldmath$n$}^T,+\mbox{\boldmath$n$}^S)$-umbilical
or $(\mbox{\boldmath$n$}^T,-\mbox{\boldmath$n$}^S)$-umbilical spacelike surface with the vanishing principal curvature in Minkowski space-time is a spacelike surface in
a lightlike hyperplane. 
Therefore we have the following proposition.
\begin{Pro}
A spacelike surface $M$ is marginally trapped with \newline $H_\ell(\bn^T,\sigma\bn^S)=0$ and totally $(\mbox{\boldmath$n$}^T,\sigma \mbox{\boldmath$n$}^S)$-umbilical
if and only if $M$ is a spacelike surface in a lightlike hyperplane, where $\sigma =+$ or $\sigma =-.$
Moreover, $M$ is strongly marginally trapped and totally $(\mbox{\boldmath$n$}^T,\pm\mbox{\boldmath$n$}^S)$-umbilical if and only if it is a spacelike plane.
\end{Pro}
\demo
By definition, we have 
\[
H_\ell (\bn^T,\sigma \bn^S)=\frac{(\kappa _1(\bn^T,\sigma \bn^S)(p)+\kappa _2(\bn^T,\sigma \bn^S)(p))}{2}.
\]
If $M$ is marginally trapped with $H_\ell(\bn^T,\sigma \bn^S)=0$ and totally $(\mbox{\boldmath$n$}^T,\sigma\mbox{\boldmath$n$}^S)$-umbilical, then $\kappa_1 (\bn^T,\sigma\bn^S)=\kappa _2(\bn^T,\sigma\bn^S)=0.$
Therefore, $M$ is a spacelike surface in a lightlike hyperplane. By the similar arguments to the above, other assertions hold.\enD
By the above proposition, we have the following corollary.
\begin{Co}\label{cor:L-flat-MT}
A spacelike surface $M$ is marginally trapped and totally \newline $(\bn^T,\bn^S)$-umbilical if and only if
up to rigid motions of $\R^4_1,$ $M$ is given by
\[
\bX_f(u_1,u_2)=(f(u_1,u_2),f(u_1,u_2),u_1,u_2)),
\]
for a smooth function $f:\R^2\lon\R.$ 
\end{Co}
\demo
By a straight forward calculation, we can show that $\bX_f$ is always spacelike.
If we consider the lightlike vector $\bv=(1,1,0,0),$ then
we have $\langle \bX_f(u_1,u_2),\bv\rangle_1= 0,$ so that $\bX_f$ is 
a spacelike surface in the lightlike hyperplane $HP(\bv,0).$
For the converse,if $M$ is a spacelike surface in a lightlike hyperplane, by a rigid motion of 
$\R^4_1$, $M$ can be included in $HP(\bv,0).$
We consider a projection $\pi :HP(\bv,0)\lon \R^2$ defined by
$\pi (x_0,x_0,x_1,x_2)=(x_1,x_2)$.
Then the fiber of $\pi$ is directed to the lightlike direction $\bv.$
Since $M$ is spacelike, $\pi|M:M\lon \R~2$ is a diffeomorphism, so that there exists a smooth
function $f:\R^2\lon\R$ such that $\bX_f(\R^2)=M.$ This completes the proof.
\enD
\begin{Rem}\rm
The above corollary indicates that there are so many complete marginally trapped surfaces in Minkowski space-time compared with minimal surfaces in Euclidean space etc.
In \cite{ChI91}, Chen and Ishikawa have shown that if $M$ is a biharmonic surface in $\R^4_1$ with flat normal connection
and marginally trapped, then up to rigid motions $M$ is given by the surface
$\bX_f(\R^2)$ in the above corollary.
We remark that if $M$ is a subset of a lightlike hyperplane, the normal connection is flat.
Related results and problems are presented in the survey article \cite{Ch09}.
\end{Rem}
\par
We deduce now the \textit{lightcone Weingarten formula:}
Since $\mbox{\boldmath$X$} _{u_i}$   
$(i=1,2$ are spacelike vectors, we have a Riemannian metric
(the {\it hyperbolic first fundamental form \/}) on $M$ defined by
$$ds^2 =g_{11}du_1^2+ g_{12}du_1du_2 + g_{22}du_2^2,$$
where $g_{ij}(u) =\langle \mbox{\boldmath$X$} _{u_i}(u ),\mbox{\boldmath$X$} _{u_j}(u)\rangle_1$ for any $u\in U.$

We also have the \textit{lightcone second fundamental form}, which is defined as the second fundamental
form associated to the normal vector field $\mbox{\boldmath$n$}^T + \mbox{\boldmath$n$} ^S $. This is given by
$$h_{ij}(\mbox{\boldmath$n$}^T,\pm\mbox{\boldmath$n$}^S )(u)=\langle -(\mbox{\boldmath$n$}^T \pm\mbox{\boldmath$n$}^S)
_{u_i}(u),\mbox{\boldmath$X$}_{u_j}(u)\rangle_1; i=1,2,$$ for any $u\in U.$
In \cite{IR07}, we have shown the following proposition.
\begin{Pro}[The lightcone Weingarten formula]
\label{p32} Under the above notations, we have the following
lightcone Weingarten formula with respect to $(\mbox{\boldmath$n$}^T,\pm\mbox{\boldmath$n$}^S)$:
\par
{\rm a)}
  $(\mbox{\boldmath$n$}^T \pm\mbox{\boldmath$n$}^S)_{u_i}
  =\langle \mbox{\boldmath$n$} ^S,\mbox{\boldmath$n$}^T_{u_i}
  \rangle_1(\mbox{\boldmath$n$}^T\pm\mbox{\boldmath$n$}^S)
  - \sum_{j=1}^2 h_i^j(\mbox{\boldmath$n$}^T,\pm\mbox{\boldmath$n$}^S )\mbox{\boldmath$X$}_{u_j}; i=1,2,$
\par
{\rm b)}
  $ \pi ^t\circ (\mbox{\boldmath$n$}^T \pm\mbox{\boldmath$n$}^S)_{u_i}
  =- \sum_{j=1}^2 h_i^j(\mbox{\boldmath$n$}^T,\pm\mbox{\boldmath$n$}^S )\mbox{\boldmath$X$}_{u_j}; i=1,2,$
\par\noindent
where $\displaystyle{\left(h_i^j(\mbox{\boldmath$n$}^T,\pm\mbox{\boldmath$n$}^S
)\right)=\left(h_{ik}(\mbox{\boldmath$n$}^T,\pm\mbox{\boldmath$n$}^S)\right)\left(g^{kj}\right)}$ and
$\displaystyle{\left( g^{kj}\right)=\left(g_{kj}\right)^{-1}}.$
\end{Pro}

The following corollary provides an explicit
expression of the lightcone curvature in terms of the Riemannian
metric and the lightcone second fundamental invariant.

\begin{Co}
\label{c33} Under the same notations as in the above proposition,
the lightcone Gauss-Kronecker curvature with respect to
$(\mbox{\boldmath$n$}^T,\pm\mbox{\boldmath$n$}^S) $ is given by
\begin{equation}
K_\ell (\mbox{\boldmath$n$}^T,\pm\mbox{\boldmath$n$}^S )=\frac{\displaystyle{{\rm
det}\left(h_{ij}(\mbox{\boldmath$n$}^T,\pm\mbox{\boldmath$n$}^S )\right)}} {\displaystyle{{\rm
det}\left(g_{\alpha \beta}\right)}}.
\end{equation}
\end{Co}

Since $\langle -(\mbox{\boldmath$n$}^T \pm\mbox{\boldmath$n$}^S )(u),\mbox{\boldmath$X$} _{u_j}(u)\rangle_1 =0,$ we have that
$$h_{ij}(\mbox{\boldmath$n$} ^T,\pm\mbox{\boldmath$n$}^S)(u)=\langle \mbox{\boldmath$n$}^T (u)\pm\mbox{\boldmath$n$}^S (u),\mbox{\boldmath$X$}
_{u_iu_j}(u)\rangle_1.$$ So the lightcone second fundamental
form at  $p_0=\mbox{\boldmath$X$} (u_0)$ depends only on the values
of the vector fields $\mbox{\boldmath$n$}^T +\mbox{\boldmath$n$}^S $ and $\mbox{\boldmath$X$} _{u_iu_j}$ at
the point $p_0$. Consequently the lightcone curvature  depends  only on
$\mbox{\boldmath$n$}^T (u_0)\pm\mbox{\boldmath$n$}^S (u_0)$, $\mbox{\boldmath$X$}_{u_i}(u_0)$  and $\mbox{\boldmath$X$}
_{u_iu_j}(u_0)$ regardless the  choice of the normal
vector fields $\mbox{\boldmath$n$}^T$ and $\mbox{\boldmath$n$}^S .$
We write as $K_\ell (\mbox{\boldmath$n$}
^T_0,\pm\mbox{\boldmath$n$}^S_0)(u_0)$ the lightcone curvature at $p_0$ with respect
to $(\mbox{\boldmath$n$} ^T_0,\pm\mbox{\boldmath$n$}^S_0)=(\mbox{\boldmath$n$}^T (u_0),\pm\mbox{\boldmath$n$}^S(u_0)).$ It thus makes sense to say
that a point $p_0$ is {\it $(\mbox{\boldmath$n$} ^T_0,\pm\mbox{\boldmath$n$}^S_0)$-umbilic} for the
lightcone $(\mbox{\boldmath$n$}^T,\pm\mbox{\boldmath$n$}^S) $-shape operator at $p_0$ just depends on
the normal vectors $(\mbox{\boldmath$n$} ^T_0,\pm\mbox{\boldmath$n$}^S_0).$
Analogously, we say that the point $p_0$ is a {\it $(\mbox{\boldmath$n$}
^T_0,\pm\mbox{\boldmath$n$}^S_0)$-parabolic point \/} of $M$ if $K_\ell (\mbox{\boldmath$n$}
^T_0,\pm\mbox{\boldmath$n$}^S_0)(u_0)=0.$ And we say that  $p_0$ is a {\it $(\mbox{\boldmath$n$}
^T_0,\pm\mbox{\boldmath$n$}^S_0)$-flat point \/} if it is  $(\mbox{\boldmath$n$} ^T_0,\pm\mbox{\boldmath$n$}^S_0)$-umbilic
 and $K_\ell(\mbox{\boldmath$n$}^T _0,\pm\mbox{\boldmath$n$}^S_0)(u_0)=0.$
 \par
 On the other hand,
 we define a normal vector field \[
 \mathfrak{K}(p)=\det S_p(\bn^T)\bn^T(u)+\det S_p(\bn^S)\bn^S(u),
 \] where $p=\bX(u).$
 We call $\mathfrak{K}$ the {\it Gaussian curvature} vector field along $M=\bX(U).$
 We also write $K(\bn^T)(p)=\det S_p(\bn^T)$ and $K(\bn^S)(p)=\det S_p(\bn^S).$
 We respective call $K(\bn^T)(p)$ and $K(\bn^S)(p)$ the {\it Gauss-Kronecker curvature} of
 $M=\bX(U)$ at $p$ with respect to $\bn^T$ and $\bn^S.$
 The {\it mean curvature} of $M=\bX(U)$ at $p$ with respect to $^bn^T$ (respectively, $\bn^S$)
 is defined to be $H(\bn^T)(p)=({\rm Trace}\, S_p(\bn^T))/2$ (respectively, $H(\bn^S)(p)=({\rm Trace}\, S_p(\bn^S))/2$)
 We also define the {\it second fundamental invariants} $h_{ij}(\bn^T)=-\langle \bn^T_{u_i},\bX_{u_j}\rangle _1$
 with respect to $\bn^T$
 and $h_{ij}(\bn^S)=-\langle \bn^S_{u_i},\bX_{u_j}\rangle _1$ with respect to $\bn^S$ respectively.
 By the standard arguments, we have the following
Weingarten formulae.
\begin{Pro} We have
\[
\pi ^t\circ\bn^T_{u_i}=-\sum _{j=1}^2 h^j_i(\bn^T)\bX_{u_j}\ {\rm and}\ \pi ^t\circ\bn^S_{u_i}=-\sum _{j=1}^2 h^j_i(\bn^S)\bX_{u_j},
\]
where $(h^j_i(\bn^T))=(h_{ik}(\bn^T))(g^{kj})$ and $(h^j_i(\bn^S))=(h_{ik}(\bn^S))(g^{kj}).$
\end{Pro}
Thus, we have the following corollary.
\begin{Co} We have
\[
K(\bn^T)=\frac{\det (h_{ij}(\bn^T))}{\det (g_{\alpha\beta})}\ {\rm and}\ K(\bn^S)=\frac{\det (h_{ij}(\bn^S))}{\det (g_{\alpha\beta})}.
\]
\end{Co}
\par
\par
We also get in this context the {\it lightcone Gauss equations} as we
shall see next.
Since $\bX (U)=M$ is a Riemannian manifold, it makes sense to consider the
{\it Christoffel symbols}:
\[
{k\brace i\ j}=\frac{1}{2}\sum_m g^{km}\left\{\frac{\partial
g_{jm}}{\partial u_i}+\frac{\partial g_{im}}{\partial u_j}
-\frac{\partial g_{ij}}{\partial u_m}\right\}.
\]
\begin{Pro} Let $\bX :U\lon \R^4_1$ be a spacelike surface. Then we have the following lightcone Gauss
equations:
\[
\bX _{u_iu_j}=\sum_k {k\brace i\ j}\bX _{u_k}-h_{ij}(\bn^T)\bn^T +h_{ij}(\bn^S)\bn^S.
\]
\end{Pro}
\demo
Since $\{\bn^T,\bn^S, \bX_{u_1},\bX_{u_2}\}$ is a 
frame of $\R^{4}_1,$
we can write
$
\bX _{u_iu_j}=\sum_k \Gamma ^k_{ij}\bx _{u_k}+\Gamma _{ij}\bn^T +\Lambda_{ij}\bn^S .
$
We now have
$$
\langle \bX _{u_iu_j},\bX _{u_{\ell}}\rangle _1 =\sum _k \Gamma^k_{ij}\langle
\bX _{u_k},\bX _{u_{\ell}}\rangle _1 =\sum_k \Gamma _{ij}^kg_{k\ell}.
$$
Since
 $
\frac{\partial g_{i\ell}}{\partial u_j}=\langle \bX _{u_iu_j},\bX
_{u_\ell}\rangle _1 +\langle \bX _{u_i},
\bX_{u_\ell u_j}\rangle _1
$
and $\bX _{u_iu_j}=\bX _{u_ju_i},$  we get
$
\Gamma ^k_{ij}=\Gamma ^k_{ji},\ \Gamma _{ij}=\Gamma _{ji},\ \Gamma
^{ij}=\Gamma ^{ji}.
$
Then by exactly the same calculation as those applied in the case of
surfaces in Euclidean space, it
follows
$
\displaystyle{\Gamma ^k_{ij}={k\brace i\ j}}.
$
\par
On the other hand,
$
-\Gamma _{ij}=\langle \bX _{u_iu_j},\bn^T\rangle _1 =h_{ij}(\bn^T).
$
Moreover
$
\Lambda_{ij}=\langle \bX _{u_iu_j},\bn^S\rangle _1 =h_{ij}(\bn^S).
$
This completes the proof.
\enD
Consider the Riemannian curvature tensor
\[
R^\ell_{ijk}=\frac{\partial }{\partial u_k}{\ell\brace i\ j}-\frac{\partial}{\partial u_j}{\ell\brace i\ k}
+\sum_m {m\brace i\ j}{\ell\brace m\ k}-\sum_m {m\brace i\ k}{\ell\brace m\ j}.
\]
Similar calculations
to those of the classical differential geometry on
surfaces in Euclidean space and the fact $\bX_{u_iu_ju_k}=\bX_{u_iu_ku_j}$
lead to the formula
\begin{eqnarray*}
R^\ell_{ijk}&{}&\!\!\!\!\!\!\!\!\! =\sum_{a} \left\{-h_{ij}(\bn^T)h_{ka}(\bn^T)+h_{ij}(\bn^S)h_{ka}(\bn^S)\right.
\\
&{}&\left.+h_{ik}(\bn^T)h_{ja}(\bn^T)-h_{ik}(\bn^S)h_ja(\bn^S)\right\}g^{a\ell}.
\end{eqnarray*}
We also consider the tensor $R_{ijk\ell}=\sum_m g_{im}R^m_{jk\ell}.$
Then we have the following proposition.
\begin{Pro}
We have
\begin{eqnarray*}
R_{ijk\ell}&{}&\!\!\!\!\!\!\!\!\! =-h_{jk}(\bn^T)h_{i\ell}(\bn^T)+h_{j\ell}(\bn^T)h_{ik}(\bn^T)\\
&{}&+h_{jk}(\bn^S)h_{i\ell}(\bn^S)-h_{j\ell}(\bn^S)h_{ik}(\bn^S).
\end{eqnarray*}
\end{Pro}
The {\it intrinsic Gauss curvature} (or, {\it sectional curvature}) of $M=\bX(U)$ at $p=\bX(u)$ is
defined by $K_I(p)=-R_{1212}(u)/\det (g_{\alpha\beta}(u)).$
 We can show that the following \lq\lq Theorema Egregium\rq\rq.
 \begin{Th} Let $K_I$ be the intrinsic Gauss curvature of $M=\bX(U).$ Then we have
 \[
 K_I(p)=-K(\bn^T)(p)+K(\bn^S)(p).
 \]
 \end{Th}
 \demo
 By Propositions 3.11, we have
 \[
 R_{1212}=\det (h_{ij}(\bn^T))-\det (h_{ij}(\bn^S)).
 \]
 Thus, by Corollary 3.9, we have
 \[
 K_I=-\frac{\det ( h_{ij} (\bn^T ))}{\det ( g_{\alpha\beta})}+\frac{\det ( h_{ij} (\bn^S))}{\det ( g_{\alpha\beta})}=-K(\bn^T)+K(\bn^S).
 \]
 This completes the proof.
 \enD
 \par
 We give the following characterization of the intrinsic flat surface (i.e., $K_I\equiv 0$) by using the
 Gaussian curvature vector $
 \mathfrak{K}(p)=K(\bn^T)(p)\bn^T(p)+K(\bn^S)(p)\bn^s(p).$
 \begin{Pro}
 For a spacelike surface $M=\bX(U)$ in $\R^4_1,$ the following conditions are equivalent{\rm:}
 \par\noindent
 {\rm 1)} $K_I(p)=0$ at $p=\bX(u).$
 \par
 \noindent
 {\rm 2)} $\mathfrak{K}(p)=\bo$ or $\mathfrak{K}(p)$ is parallel to the lightlike vector
 $\bn^T(p)+\bn^S(p).$
 \end{Pro}
 \demo
 Since $K_I=-K(\bn^T)+K(\bn^S),$ 
 $K_I(p)=0$ if and only if $K(\bn^T)(p)=K(\bn^T)(p).$
 Therefore, we have 
\[\mathfrak{K}(p)=K(\bn^T)(p)\bn^T(p)+K(\bn^S)(p)\bn^S(p)
 = K(\bn^T)(\bn^T(p)+\bn^S(p)).
 \]
 \par
 For the converse, if $\mathfrak{K}(p)=\lambda (\bn^T(p)+\bn^S(p)),$
 then $K(\bn^T)(p)=\lambda =K(\bn^S)(p),$ so that $K_I(p)=0.$
 \enD
 We also say that $M=\bX(U)$ is an {\it extrinsic flat surface} if $\mathfrak{K}\equiv \bo.$

 \section{Normalized lightcone curvatures}

 Given a spacelike embedding $\mbox{\boldmath$X$}:U\rightarrow  \mathbb{R}^{4}_1$  from an open
subset $U\subset \mathbb{R}^{2},$ and a point $p=\mbox{\boldmath$X$}(u)$,  consider a
future directed unit timelike normal section $\mbox{\boldmath$n$} ^T(u)\in N_p(M)$
and the corresponding spacelike unit normal section $\mbox{\boldmath$n$}^S(u)\in
N_p(M)$ constructed in the previous section. Given another
future directed unit timelike normal section $\overline{\mbox{\boldmath$n$}}^T(u),$ we
have
$\widetilde{(\mbox{\boldmath$n$}^T\pm\mbox{\boldmath$n$}^S)}(u)=\widetilde{(\overline{\mbox{\boldmath$n$}}^T\pm\overline{\mbox{\boldmath$n$}}^S)}(u)\in
S^{2}_+,$ so we have a well-defined {\it lightcone Gauss map}
of $M=\mbox{\boldmath$X$}(U)$ as
$$\begin{array}{cccl} \widetilde{\mathbb L}^\pm : &  U & \longrightarrow & S^{2}_+ \\
 & u& \longmapsto & \widetilde{(\mbox{\boldmath$n$}^T\pm\mbox{\boldmath$n$}^S)}(u).
\end{array}$$
This induces a linear mapping $d\widetilde{\mathbb L}^\pm_p:T_pM\rightarrow
T_p \mathbb{R}^{4}_1$ under the identification of $U$ and $M,$ where
$p=\mbox{\boldmath$X$}(u).$
The  following proposition has been shown in \cite{IR07}.

\begin{Pro}
With the above notation, we have the following normalized lightcone
Weingarten formula:
\begin{equation}
\pi^t\circ\widetilde{\mathbb L}^\pm_{u_i}=-\sum_{j=1}^{2}
\frac{1}{\ell ^\pm_0(u)}h^j_i(\mbox{\boldmath$n$}^T,\pm\mbox{\boldmath$n$}^S)\mbox{\boldmath$X$}_{u_j}, 
\end{equation}
where ${\mathbb L}^\pm(u)=(\ell ^\pm_0(u),\ell ^\pm_1(u),\ell ^\pm_2(u) ,\ell ^\pm_3(u)).$
\end{Pro}
\par
We now define the {\it normalized lightcone second fundamental invariant} by
\[
\widetilde{h}[\pm]_{ij}(u)=\frac{1}{\ell _0(u)}h_{ij}(\bn^T,\pm\bn^S)(u).
\]
Since $h_{ij}(\bn^T,\pm\bn^S)(u)=\langle (\mbox{\boldmath$n$}^T \pm\mbox{\boldmath$n$}^S)
(u),\mbox{\boldmath$X$}_{u_iu_j}(u)\rangle_1,$
we have
\begin{eqnarray*}
&{}& \widetilde{h}[\pm]_{ij}(u)=\frac{1}{\ell _0(u)}\langle (\mbox{\boldmath$n$}^T \pm\mbox{\boldmath$n$}^S)
_{u_i}(u),\mbox{\boldmath$X$}_{u_j}(u)\rangle _1 \\
&{}& =\langle \widetilde{\mathbb{L}}^\pm(u),\mbox{\boldmath$X$}_{u_iu_j}(u)\rangle _1  =\langle -\widetilde{\mathbb{L}}^\pm_{u_i}(u),\mbox{\boldmath$X$}_{u_j}(u)\rangle _1.
\end{eqnarray*}

We need the following detailed formula.
\begin{Lem} 
With the above notations, we have 
\[
\widetilde{\mathbb{L}}^\pm_{u_i}=\left(\frac{-\ell _{0u_i}}{\ell _0^2}+\frac{\langle \bn^S,\bn^T_{u_i}\rangle _1}{\ell _0}\right)(\bn^T\pm \bn^S)-\sum_{j=1}^{2}
\widetilde{h}[\pm]^j_i\mbox{\boldmath$X$}_{u_j},
\]
where $(\widetilde{h}[\pm]^j_i)=(\widetilde{h}[\pm]_{ik})(g^{kj}).$
\end{Lem}
\demo
\par
In the proof of Proposition 4.1, we have 
$
\ell _0\widetilde{\mathbb{L}}^\pm_{u_i}=(\bn^T\pm \bn^S)_{u_i}-\ell _{0u_i}\widetilde{\mathbb{L}}^\pm.
$
Thus the above formula directly follows form the assertion a) of Proposition 3.3 and the relation $\ell _0\widetilde{\mathbb{L}}^\pm=\bn^T\pm \bn^S.$
\enD
We call the linear transformation $\widetilde{S}^\pm_p=-\pi^t\circ
d\widetilde{\mathbb L}^\pm_p$ the {\it normalized lightcone shape
operator} of $M=\mbox{\boldmath$X$} (U)$ at $p.$ The {\it normalized lightcone
Gauss-Kronecker curvature} of $M=\mbox{\boldmath$X$} (U)$ is defined to be
$\widetilde{K}^\pm_\ell (u)={\rm det}\, \widetilde{S}^\pm_p.$ We say that
$p=\mbox{\boldmath$X$}(u)$ is a {\it lightlike parabolic point \/} if
 $\widetilde{K}^\pm_\ell (u)=0.$
We also define the {\it normalized lightcone mean curvature} of $M=\mbox{\boldmath$X$} (U)$
by $\widetilde{H}^\pm_\ell (p)={\rm Trace }\, \widetilde{S}^\pm_p/2.$
The eigenvalues $\{ \widetilde{\kappa}^\pm_i(p)\}_{i=1}^{2}$ of
$\widetilde{S}^\pm_p$ are called {\it normalized lightcone principal
curvatures}. It follows from the above formula that
$\widetilde{\kappa}^\pm _i(p)=(1/\ell _0)\kappa _i(\mbox{\boldmath$n$}^T,\pm\mbox{\boldmath$n$}^S)(p).$
Clearly, the eigenvectors of $\widetilde{S}^\pm_p$ coincide with the
lightcone principal directions with respect to $(\mbox{\boldmath$n$}^T,\pm\mbox{\boldmath$n$}^S) $, for
any future directed frame
 $(\mbox{\boldmath$n$}^T,\pm\mbox{\boldmath$n$}^S) $ on $M$, therefore, we can refer to the
 $(\mbox{\boldmath$n$}^T,\pm\mbox{\boldmath$n$}^S) $-lightcone principal configuration, simply as the
 {\it lightcone principal configuration} on $M$.  The $(\mbox{\boldmath$n$}^T,\pm\mbox{\boldmath$n$}^S)
 $-umbilics  shall be called {\it lightlike umbilics}.
We say that $M=\mbox{\boldmath$X$} (U)$ is {\it totally lightlike umbilical} if all
points on $M$ are lightlike umbilic. The point $p$ is
called a {\it lightlike flat point} if $p$ is both lightlike
umbilic and parabolic.
 The spacelike submanifold $M=\mbox{\boldmath$X$}(U)$ is called  {\it lightlike flat}
 provided every point of $M$ is lightlike flat. As observed in the previous section,
 the lightcone principal configuration is preserved by Lorentz transformations,
 although the normalized lightcone principal curvatures are not.
 We remark that the normalized principal curvatures are invariant under the $SO(3)$-action and parallel translations, where
 $SO(3)$ is the canonical subgroup of $SO_0(1,3).$
 \par
 By Proposition 4.1, we have the following relations:
 \[
 \widetilde{K}^\pm_\ell(p)=\left(\frac{1}{\ell^\pm_0(p)}\right)^2K^\pm_\ell (\mbox{\boldmath$n$}^T,\pm \mbox{\boldmath$n$}^S)(p),\quad
 \widetilde{H}^\pm_\ell (p)=\frac{1}{\ell^\pm_0(p)}H_\ell (\mbox{\boldmath$n$}^T,\pm \mbox{\boldmath$n$}^S)(p).
 \]
 Therefore, we have the following proposition.
 \begin{Pro} Let $M=\mbox{\boldmath$X$}(U)$ be a spacelike surface in $\mathbb{R}^4_1.$ Then
 \par
 {\rm 1)} $\widetilde{K}^\pm_\ell(p)=0$ if and only if $K^\pm_\ell (\mbox{\boldmath$n$}^T,\pm \mbox{\boldmath$n$}^S)(p)=0$ for any future directed normal frame $(\mbox{\boldmath$n$}^T,\mbox{\boldmath$n$}^S).$
 \par
 {\rm 2)} $\widetilde{H}^\pm_\ell (p)=0$ if and only if $H_\ell (\mbox{\boldmath$n$}^T,\pm \mbox{\boldmath$n$}^S)(p)=0$ for any future directed normal frame $(\mbox{\boldmath$n$}^T,\mbox{\boldmath$n$}^S).$
\end{Pro}
We have the following corollary.
\begin{Co} Let $M=\mbox{\boldmath$X$}(U)$ be a spacelike surface in $\mathbb{R}^4_1.$
Then $M$ is a marginally trapped surface if and only if $\widetilde{H}^+_\ell \equiv 0$ or $\widetilde{H}^-_\ell \equiv 0$.
Moreover, $M=\bX(U)$ is strongly marginally trapped surface if and only if $\widetilde{H}^+_\ell\equiv \widetilde{H}^-_\ell \equiv 0$.
\end{Co} 
\par
On the other hand, we say that $(u_1,u_2)\in U$ is an {\it isothermal parameter} of $M=\bX(U)$ if
$g_{11}=g_{22}$ and $g_{12}=0.$
For any two dimensional Riemannian manifold, there exists an isothermal parameter at any point \cite{Cher}.
\begin{Pro}
Let $\bX:U\lon \R^4_1$ be a spacelike surface with an isothermal parameter $(u_1,u_2)\in U.$
Then we have
\[
\bX_{u_1u_1}+\bX_{u_2u_2}=g_{11}\left(-\ell ^+_0(\widetilde{H}^+_\ell+\widetilde{H}^-_\ell)\bn^T
+\ell^-_0(\widetilde{H}^+_\ell-\widetilde{H}^-_\ell)\bn^S\right).
\]
\end{Pro}
\demo
Since $(u_1,u_2)$ is an isothermal parameter, we have
\[
\langle \bX_{u_1},\bX_{u_1}\rangle _1 = \langle \bX_{u_2},\bX_{u_2}\rangle _1,\ \langle \bX_{u_1},\bX_{u_2}\rangle _1=0.
\]
It follows that
\[
\langle \bX_{u_1u_1},\bX_{u_1}\rangle _1 = \langle \bX_{u_1u_2},\bX_{u_2}\rangle _1= -\langle \bX_{u_1},\bX_{u_2u_2}\rangle _1.
\]
Therefore, we have $\langle \bX_{u_1u_1}+\bX_{u_2u_2},\bX_{u_1}\rangle _1=0.$
By similar arguments to the above, we have
$\langle \bX_{u_1u_1}+\bX_{u_2u_2},\bX_{u_2}\rangle _1=0.$
Thus the vector $\bX_{u_1u_1}+\bX_{u_2u_2}$ is normal to $M=\bX(U),$ so that there exist $\lambda, \mu \in \R$ 
such that $\bX_{u_1u_1}+\bX_{u_2u_2}=\lambda\bn^T+\mu\bn^S.$
Then we have
\begin{eqnarray*}
-\lambda &{}&\!\!\!\!\!\!\! =\langle \bX_{u_1u_1}+\bX_{u_2u_2},\bn^T\rangle _1 =h_{11}(\bn^T)+h_{22}(\bn^T) \\
\mu &{}&\!\!\!\!\!\!\! = \langle \bX_{u_1u_1}+\bX_{u_2u_2},\bn^S\rangle _1 =h_{11}(\bn^S)+h_{22}(\bn^S).
\end{eqnarray*}
Since $\ell^\pm_0\widetilde{h}[\pm]_{ij}=h_{ij}(\bn^T,\pm\bn^S)=h_{ij}(\bn^T)\pm h_{ij}(\bn^S),$ we get
\[
\ell^\pm_0(\widetilde{h}[\pm]_{22}+\widetilde{h}[\pm]_{11})
=h_{22}(\bn^T)+h_{11}(\bn^T)\pm (h_{22}(\bn^S)+h_{11}(\bn^S)).
\]
Therefore we have
\begin{eqnarray*}
2(h_{11}(\bn^T)+h_{22}(\bn^T))&{}&\!\!\!\!\!\!\!=\ell ^+_0(\widetilde{h}[+]_{22}+\widetilde{h}[+]_{11}+\widetilde{h}[-]_{22}+\widetilde{h}[-]_{11})\\
&{}&\!\!\!\!\!\!\!=2\ell^+_0g_{11}(\widetilde{H}^+_\ell +\widetilde{H}^-_\ell) \\
2(h_{11}(\bn^S)+h_{22}(\bn^S))&{}&\!\!\!\!\!\!\!=\ell ^-_0(\widetilde{h}[+]_{22}+\widetilde{h}[+]_{11}-\widetilde{h}[-]_{22}-\widetilde{h}[-]_{11}) \\
&{}&\!\!\!\!\!\!\!=2\ell^-_0g_{11}(\widetilde{H}^+_\ell -\widetilde{H}^-_\ell)
\end{eqnarray*}
This means that
\[
\bX_{u_1u_1}+\bX_{u_2u_2}=-\ell^+_0g_{11}(\widetilde{H}^+_\ell +\widetilde{H}^-_\ell)\bn^T+\ell^-_0g_{11}(\widetilde{H}^+_\ell -\widetilde{H}^-_\ell)\bn^S.
\]
This completes the proof.
\enD
\begin{Co}
Let $\bX:U\lon \R^4_1$ be a spacelike surface with an isothermal parameter $(u_1,u_2)\in U.$
Then $\bX_{u_1u_1}+\bX_{u_2u_2}=\bo$ if and only if $M=\bX(U)$ is a strongly marginally trapped
surface.
\end{Co}
\demo
By the above proposition, $\bX_{u_1u_1}+\bX_{u_2u_2}=\bo$ if and only if
$\widetilde{H}^+_\ell +\widetilde{H}^-_\ell=\widetilde{H}^+_\ell -\widetilde{H}^-_\ell=0.$
The last conditions mean that $\widetilde{H}^+_\ell =\widetilde{H}^-_\ell=0.$
Since $\ell ^\pm _o\widetilde{H}^\pm_\ell=H(\bn^T,\pm\bn^S)=H(\bn^T)\pm H(\bn^S),$
the last condition is equivalent to the condition
$H(\bn^T)=H(\bn^S)=0$ which means $\mathfrak{H}=0.$
\enD

%%%%%%%%%%%%%%%%%%%%%%%%%%%%%%%%%%
\section{Variation formula for marginally trapped surfaces}

In this section, we prove the first and second variation formula of area
(Proposition \ref{Pro:1st-var} and Theorem \ref{thm:SVF} ).
Let $U$ be a bounded compact domain in $\mathbb{R}^2$. 
The area $A(\bX)$ of a spacelike embedding 
$\bX : U \rightarrow \mathbb{R}^4_1$
is given by 
\[
  A(\bX) 
  = \int_{U} dM,\qquad 
  \left(dM= \sqrt{g_{11}g_{22} - g_{22}^2}\, du_1du_2\right).
\]
We denote by $\{\bn^S,\,\bn^T\}$ an orthonormal frame of the normal bundle of $\bX$, 
where $\bn^S$ and $\bn^T$ are spacelike and timelike, respectively. 
Then, if we set the lightlike vectors $\bl^{\pm}$ as 
$\bl^{\pm} = \bn^T \pm \bn^S$, 
then the pair $\{\bl^{+},\, \bl^{-}\}$ defines a basis of $NU$, 
and satisfies the following relations
\[
  \inner{\bl^{+}}{\bl^{+}}=
  \inner{\bl^{-}}{\bl^{-}}=0,\qquad
  \inner{\bl^{+}}{\bl^{-}}=-2.
\]

%%%%%%%%%%%%%%%%%%%%%%%%%%%%%%%%%%
\subsection{First variation formula}

Let $\bX^{\ve}$ be a smooth variation of 
a spacelike embedding $\bX : U\rightarrow \mathbb{R}^4_1$,
where $\ve$ denotes the variation parameter.
That is, 
$\bX^{\ve} : U \rightarrow \mathbb{R}^4_1$
is a smooth spacelike surface for each $\ve$,
and it satisfies $\bX^{0}=\bX$.

We assume that the variation vector field $V$ of the variation $\bX^{\ve}$ is given by
\begin{equation}\label{eq:VVF}
  V^{\pm}= \Pdeo \bX^{\ve} = \alpha \bl^{\pm}
\end{equation}
where $\alpha$ is a smooth function on $U$.

\begin{Pro}[First variation formula]\label{Pro:1st-var}
For a spacelike surface
$\bX:\Sigma\rightarrow \mathbb{R}^4_1$,
let $\bX^{\varepsilon}$ be a variation
whose variation vector field $V$ of $\bX^{\varepsilon}$ is given by \eqref{eq:VVF}.
Then the first variation of area is given by 
\[
  \left. \frac{d}{d \ve} \right|_{\ve=0} A(\bX^{\ve})
  = -2 \int_{U} \alpha H_{\ell}(\bn^{T}, \pm \bn^{S}) dM.
\]
\end{Pro}

\demo
For the sake of simplicity, we will give the proof in the case of 
$V=V^{+}=\alpha \bl^{+}$ here. 
Similar proof is valid for the case of $V^{-}$.
Since the area $A(\bX^{\varepsilon})$ of $\bX^{\varepsilon}$ is 
\[
  A(\bX^{\varepsilon})
  = \int_{U} \sqrt{g_{11}^{\varepsilon}g_{22}^{\varepsilon}-(g_{22}^{\varepsilon})^2}\,du_1du_2,
\]
we have
\begin{align*}
  \frac{d}{d \varepsilon} A(\bX^{\varepsilon})
  &=\int_{U} \Pde 
       \left( \sqrt{g_{11}^{\varepsilon}g_{22}^{\varepsilon}-(g_{22}^{\varepsilon})^2} \right)du_1du_2 \\
  &= \int_{U}
       \frac{\de \left(g_{11}^{\varepsilon}g_{22}^{\varepsilon}-(g_{22}^{\varepsilon})^2\right)
       }{ 2 \left(g_{11}^{\varepsilon}g_{22}^{\varepsilon}-(g_{22}^{\varepsilon})^2\right) } dM,
\end{align*}
where $\partial_{\ve} = \partial/\partial \ve$. Since
\[
  \Pde 
       \left(g_{11}^{\varepsilon}g_{22}^{\varepsilon}-(g_{22}^{\varepsilon})^2\right)
  = \left( \Pde g_{11}^{\varepsilon} \right) g_{22}^{\varepsilon}
     + g_{11}^{\varepsilon} \left( \Pde g_{22}^{\varepsilon} \right)
     - 2g_{12}^{\varepsilon} \left( \Pde g_{12}^{\varepsilon} \right),
\]
and
\begin{align*}
   \Pdeo g_{ij}^{\varepsilon}
  &= \Pdeo 
        \inner{ \bX^{\varepsilon}_{u_i} }{ \bX^{\varepsilon}_{u_j} }\\
  &= \inner{ V_{u_i} }{ \bX_{u_j} } 
        + \inner{ \bX_{u_i} }{ V_{u_j} }\\
  &= \inner{ \alpha \bl^{+}_{u_i} }{ \bX_{u_j} } 
        + \inner{ \bX_{u_i} }{ \alpha \bl^{+}_{u_j} }
  = - 2\alpha \hpij,
\end{align*}
we obtain
\begin{align*}
  \frac{\deo \left(g_{11}^{\varepsilon}g_{22}^{\varepsilon}-(g_{22}^{\varepsilon})^2\right)
        }{ 2 \left(g_{11}g_{22}-g_{22}^2\right) }
  = -\alpha \sum_{i,j} g^{ij} \hpij 
  =  -2\alpha \Hp.
\end{align*} 
This completes the proof.
\enD

%%%%%%%%%%%%%%%%%%%%%%%%%%%%%%%%%%
\subsection{Second variation formula}

Now, we prove the second variation formula for area.

\begin{Th}[Second variation formula]\label{thm:SVF}
Let $\bX:\Sigma\rightarrow \mathbb{R}^4_1$ be a marginally trapped surface,
that is, $\bX$ satisfies $H_{\ell}(\bn^{T},\bn^{S})=0$
(resp.\ $H_{\ell}(\bn^{T},-\bn^{S})=0$).
If $\bX^{\varepsilon}$ is a variation
whose variation vector field $V$ of $\bX^{\varepsilon}$ is given by 
\[
  V^{+}= \Pdeo \bX^{\ve} = \alpha \bl^{+}\qquad
  \left(\text{resp. } V^{-}= \Pdeo \bX^{\ve} = \alpha \bl^{-} \right),
\]
then the second variation of area satisfies
\begin{align*}
  \left. \frac{d^2}{d \varepsilon^2} \right|_{\varepsilon=0} A(\bX^{\varepsilon})
  &= 2 \int_{U} \alpha^2 K_{\ell}(\bn^{T},\bn^{S}) dM\\ 
  &{} \left(\text{resp. } 
  \left. \frac{d^2}{d \varepsilon^2} \right|_{\varepsilon=0} A(\bX^{\varepsilon})
  = 2 \int_{U} \alpha^2 K_{\ell}(\bn^{T},-\bn^{S}) dM \right).
\end{align*}
\end{Th}

For the proof of Theorem \ref{thm:SVF}, we need the following Lemma \ref{lem:W} and Lemma \ref{lem:mean-variation}.
The Lemma \ref{lem:W} is obtained by the Cayley-Hamilton theorem.

\begin{Lem}\label{lem:W}
For spacelike surfaces in $\mathbb{R}^4_1$, we have the following{\rm :}
\begin{itemize}
\item[$(1)$] $\sum_{k}\Wpm{i}{k} \Wpm{k}{j} = 2\Hpm \Wpm{i}{j}$ \newline $-\Kpm \delta^i_j$,
\item[$(2)$] $\inner{\bl^{\pm}_{u_i}}{\bl^{\pm}_{u_j}} = 2\Hpm \hpmij - \Kpm g_{ij}$.
\end{itemize}
\end{Lem}

\begin{Lem}\label{lem:mean-variation}
For a variation vector field $V=\alpha \bl^{+}$, we have
\begin{eqnarray}\label{eq:mean-variation}
  \Pdeo \Hp &{}&\!\!\!\!\!\!\!\!\!\!= \alpha\left(2 \Hp^2-\Kp\right) \nonumber \\ 
   &{}&\!\!\! -\frac{1}{2}\inner{\Pdeo \bl^{+}}{\bl^{-}} \Hp.
\end{eqnarray}
\end{Lem}

\demo
We have 
\begin{align*}
  2\Pdeo \Hp 
  &= \Pdeo \left( g^{ij} \hpij \right)\\ 
  &= \left(\Pdeo g^{ij}\right) \hpij + g^{ij}\left(\Pdeo \hpij\right)\\  
  &= -g^{ib} \left(\Pdeo g_{ab}\right) g^{aj} \hpij \\
  &{}\hspace{20mm}+ g^{ij} \Pdeo \inner{\bl^{+}}{\bX_{u_iu_j}}\\
  &= -g^{ib} g^{aj} \hpij \left(\Pdeo g_{ab}\right) \\
  & {}\hspace{10mm}    + g^{ij}\left(\inner{\Pdeo \bl^{+}}{\bX_{u_iu_j}} + \inner{\bl^{+}}{V_{u_iu_j}}\right)\\
  &= P_1+ P_2 + P_3,
\end{align*}
where 
\[
  P_1=-g^{ib} g^{aj} \hpij \left(\Pdeo \inner{\bX_{u_a}}{\bX_{u_b}}\right),\
  P_2=g^{ij}\inner{\Pdeo \bl^{+}}{\bX_{u_iu_j}},
\]
and $P_3= g^{ij}\inner{\bl^{+}}{V_{u_iu_j}}$.
With respect to $P_1$, it follows that
\begin{align*}
  P_1
  &= g^{ib} g^{aj} \hpij \left(\inner{V_{u_a}}{\bX_{u_b}}+\inner{\bX_{u_a}}{V_{u_b}}\right)\\
  &= 2g^{ib} g^{aj} \hpij \inner{V}{\bX_{u_au_b}}\\
  &= 2\alpha g^{aj} g^{ib} \hpij \hpab 
  = 2\alpha \Wpm{j}{b} \Wpm{b}{j}.
\end{align*}
By $(1)$ in Lemma \ref{lem:W}, we have
$\Wpm{j}{b} \Wpm{b}{j} = 4\Hp^2-2\Kp$, 
and hence  
\begin{equation}\label{eq:first-part}
  P_1 = 4\alpha \left( 2\Hp^2-\Kp \right)
\end{equation}
holds.
For $P_3$, we have
\begin{align*}
  P_3
  &= g^{ij}\inner{\bl^{+}}{ \alpha_{u_iu_j}\bl^{+} + 2\alpha_{u_i}\bl^{+}_{u_j} + \alpha \bl^{+}_{u_iu_j} }\\
  &= \alpha g^{ij} \inner{\bl^{+}}{ \bl^{+}_{u_iu_j} }
  = - \alpha g^{ij} \inner{\bl^{+}_{u_i}}{ \bl^{+}_{u_j} }.
\end{align*}
Then, by $(2)$ in Lemma \ref{lem:W}, it follows that
\begin{equation} \label{eq:third-part}
\begin{split}
  P_3
  &= - \alpha g^{ij} \left( 2\Hpm \hpmij - \Kpm g_{ij}\right)\\
  &= - 2\alpha \left(2\Hp^2 - \Kp\right).
\end{split}
\end{equation}
Finally, for $P_2$, we shall calculate $\deo \bl^{+}$.
If we put $\deo \bl^{+}= p^j \bX_{u_j} + \Phi \bl^{+}$,
we have $\inner{\deo \bl^{+}}{\bX_{u_k}} = \inner{p^j \bX_{u_j}}{\bX_{u_k}} = p^j g_{jk}$.
Hence 
\[
  g^{kl} \inner{\Pdeo \bl^{+} }{\bX_{u_k}} 
  = g^{kl} p^j g_{jk}
  = p^j \delta^l_j
  = p^l
\]
holds.
On the other hand, we have
\[
  \inner{\Pdeo \bl^{+} }{\bX_{u_k}} 
  = -\inner{\bl^{+} }{V_{u_k}}
  = -\inner{\bl^{+} }{\alpha_{u_k}\bl^{+}+\alpha \bl^{+}_{u_k}}
  = 0,
\]
and thus we have $p^{j}=0$.
Since $\inner{\deo \bl^{+}}{\bl^{-}}=-2\Phi$, 
we have
\[
  \Pdeo \bl^{+} = -\frac{1}{2} \inner{\Pdeo \bl^{+}}{\bl^{-}}\bl^{+}.
\]
Thus we get
\begin{equation}\label{eq:second-part}
  P_2  
  = -\frac{1}{2} \inner{\Pdeo \bl^{+}}{\bl^{-}} g^{ij} \hpij
  = - \inner{\Pdeo \bl^{+}}{\bl^{-}} \Hp.
\end{equation}
By \eqref{eq:first-part}, \eqref{eq:third-part} and \eqref{eq:second-part},
we obtain \eqref{eq:mean-variation}.
\enD

\noindent
{\it Proof of Theorem \ref{thm:SVF}}.\,
%\demo
We shall give the proof in the case of $\Hp=0$. 
A similar proof is valid for the case of $H_{\ell}(\bn^{T},-\bn^{S})=0$.
In this situation, the variation vector field is given by $V=V^{+}=\alpha \bl^{+}$.
By Proposition \ref{Pro:1st-var},
the first variation formula of the area is given by
\[
  \frac{d^2}{d \varepsilon^2} A(\bX^{\varepsilon})
  = \int_{U} \Hp \inner{\Pde \bX^{\varepsilon}}{\bl^{-}} dM.
\]
Since $\Hp=0$ holds at $\varepsilon=0$,
the second variation of area is given by
\begin{align*}
  \left. \frac{d^2}{d \varepsilon^2} \right|_{\varepsilon=0} A(\bX^{\varepsilon})
  &= \int_{U}\left\{ \left( \Pde \Hp \right) 
       \inner{\Pde \bX^{\varepsilon}}{\bl^{-}} dM\right.\\
  &{}\hspace{10mm} \left.\left.+ \Hp \Pde \left(\inner{\Pde \bX^{\varepsilon}}{\bl^{-}} 
           dM \right) \right\} \right|_{\varepsilon=0}\\
  &= \int_{U} \left( \Pdeo \Hp \right)
        \inner{V}{\bl^{-}} dM.
\end{align*}
Then, by Lemma \ref{lem:mean-variation}, we have
\[
  \left. \frac{d}{d \varepsilon} \right|_{\varepsilon=0} \Hp 
  = -\alpha \Kp.
\]
Since $\inner{V}{\bl^{-}} = -2\alpha$, we get the proof.
\enD

\begin{Rem}
{\rm
It should be remarked that the quantity
    $\de \Hp$ which we calculated in
    Lemma \ref{lem:mean-variation} is also calculated
    by Andersson and Metzger \cite[Lemma 4.1]{AM}
    in the case of a general Lorentzian $4$-manifold.

}
\end{Rem}

%%%%%%%%%%%%%%%%%%%%%%%%%%%%%%%%%%
\section{Marginally trapped graphs and the Bernstein-type problem}
\label{sec:graph}

In this section, we consider a system of partial differential equations for marginally trapped surfaces,
which are given in the graph form, 
satisfy the equation.
Then we have the following theorem.

\begin{Th}\label{thm:graph}
Let $U$ be a domain of $\mathbb{R}^2$, and $f,\,g$ be smooth functions on $U$.
Consider a graph immersion $\bX : U \rightarrow \mathbb{R}^4_1$ 
given by $$
\bX(u_1,u_2)=(f(u_1,u_2),g(u_1,u_2),u_1,u_2).
$$
If we set functions $\phi_1,~\phi_2$ and $\Delta$ as
\begin{align}
  \phi_1 &= (1-f_{u_2}^2+g_{u_2}^2) f_{u_1u_1} -2 (-f_{u_1} f_{u_2}+g_{u_1} g_{u_2}) f_{u_1u_2} 
  +(1-f_{u_1}^2+g_{u_1}^2) f_{u_2u_2},
  \label{eq:alpha}\\
  \phi_2 &= (1-f_{u_2}^2+g_{u_2}^2) g_{u_1u_1} -2 (-f_{u_1} f_{u_2}+g_{u_1} g_{u_2}) g_{u_1u_2} 
  +(1-f_{u_1}^2+g_{u_1}^2) g_{u_2u_2},
  \label{eq:beta}\\
  \Delta&=(f_{u_1} g_{u_1}+f_{u_2} g_{u_2})^2- (1+g_{u_1}^2+g_{u_2}^2)(-1+f_{u_1}^2+f_{u_2}^2),
  \label{eq:Delta}
\end{align}
then, $\bX$ is strongly marginally trapped if and only if
\begin{equation}\label{eq:SMTS-eq}
  \Delta>0 \qquad \text{and} \qquad
  \phi_1=\phi_2=0
\end{equation}
hold. 
Moreover, $\bX$ is marginally trapped if and only if
the functions $f,\,g$ and $\Delta$ satisfy 
\begin{equation}\label{eq:MTS-eq}
  \Delta>0 \quad \text{and} \quad
  (1+g_{u_1}^2+g_{u_2}^2)\phi_1^2 - 2(f_{u_1} g_{u_1}+f_{u_2} g_{u_2})\phi_1 \phi_2 
  + (-1+f_{u_1}^2+f_{u_2}^2)\phi_2^2=0.
\end{equation}
\end{Th}

For the proof of this theorem, we calculate the mean curvature vector 
of a graph immersion.

\begin{Lem}\label{lem:g-MCV}
For a spacelike surface $\bX : U\rightarrow \mathbb{R}^4_1$ defined by
$\bX(u_1,u_2)=(f(u_1,u_2),g(u_1,u_2),u_1,u_2)$, 
let $\phi_1,~\phi_2$ and $\Delta$ be functions defined as in
\eqref{eq:alpha}, \eqref{eq:beta} and \eqref{eq:Delta}, respectively.
Setting $\vect{v}_1$ and $\vect{v}_2$ as 
\[
  \vect{v}_1=(1,\,0,\,f_u,\,f_v),\qquad
  \vect{v}_2=(0,\,1,\,-g_u,\,-g_v),
\]
and $\tau_{ij}\, (i,j=1,2)$ as
\[
  \tau_{11} = 1+g_{u_1}^2+g_{u_2}^2,\qquad
  \tau_{12} = \tau_{21} = f_{u_1} g_{u_1}+f_{u_2} g_{u_2},\qquad
  \tau_{22} = -1+f_{u_1}^2+f_{u_2}^2,
\]
then the mean curvature vector $\mathfrak{H}$ of $\bX$ is given by
\begin{equation}\label{eq:g-MCV}
  \mathfrak{H}=\frac{\left(\tau_{11}\phi_1 - \tau_{12}\phi_2\right) \vect{v}_1
  +\left(\tau_{12}\phi_1 - \tau_{22}\phi_2\right) \vect{v}_2}{2\Delta^2}.
\end{equation}
\end{Lem}

\demo
The mean curvature vector $\mathfrak{H}$ is given by
\begin{equation}\label{eq:MCV}
  \mathfrak{H}=\frac{g_{22} I\!I( \bX _{u_1}, \bX _{u_1}) -2g_{12} I\!I( \bX _{u_1}, \bX _{u_2})
  +g_{11} I\!I( \bX _{u_2}, \bX _{u_2})}{2(g_{11}g_{22}-g_{12}^2)},
\end{equation}
where 
\begin{equation}\label{eq:2ndFF}
  I\!I( \bX _{u_i}, \bX _{u_j})= \bX _{u_iu_j} - \sum_{k} {k\brace i\ j} \bX _{u_k},
\end{equation}
for $i,j =1,2$, and ${k\brace i\ j}$ 
are the Christoffel symbols of the induced metric $g= g_{ij}du_idu_j$.
The coefficients $g_{ij}$ of the metric $g$ are calculated as
\begin{align*}
  g_{11}&=\inner{\bX_{u_1}}{\bX_{u_1}}=1-f_{u_1}^2+g_{u_1}^2,\\
  g_{12}&=\inner{\bX_{u_1}}{\bX_{u_2}}=-f_{u_1} f_{u_2}+g_{u_1} g_{u_2},\\
  g_{22}&=\inner{\bX_{u_2}}{\bX_{u_2}}=1-f_{u_2}^2+g_{u_2}^2,
\end{align*}
and hence we have
\begin{equation} \label{eq:metric}
\begin{split}
  g_{11}g_{22}-g_{12}^2
  &=(1-f_{u_1}^2+g_{u_1}^2)(1-f_{u_2}^2+g_{u_2}^2)-(-f_{u_1} f_{u_2}+g_{u_1} g_{u_2})^2\\
  &=\tau_{12}^2- \tau_{11} \tau_{22}
  =\Delta.
\end{split}
\end{equation}
Substituting these into \eqref{eq:MCV} and \eqref{eq:2ndFF},
we obtain \eqref{eq:g-MCV}.
\enD

\noindent
{\it Proof of Theorem \ref{thm:graph}}.\,
%\demo
By \eqref{eq:metric}, we have that
the graph immersion $\bX$ is spacelike if and only if $\Delta>0$.
Moreover, by Lemma \ref{lem:g-MCV} 
and the linear independentness of $\vect{v}_1$ and $\vect{v}_2$, 
it follows that the spacelike graph immersion $\bX$ is strongly marginally trapped 
if and only if 
\begin{equation}\label{eq:proof-SMT}
  \tau_{11}\phi_1 - \tau_{12}\phi_2=0,\qquad
  \tau_{12}\phi_1 - \tau_{22}\phi_2=0
\end{equation}
hold.
Since $\Delta>0$ and by \eqref{eq:metric}, we have that \eqref{eq:proof-SMT} 
is equivalent to $\phi_1=\phi_2=0$.
This completes the proof of \eqref{eq:SMTS-eq}.

With respect to \eqref{eq:MTS-eq},
by \eqref{eq:g-MCV} in Lemma \ref{lem:g-MCV}, 
we have that
the spacelike graph immersion $\bX$ is strongly marginally trapped 
if and only if 
\begin{equation}\label{eq:proof-MT}
\langle \ba,\bb\rangle_1=0,
\end{equation}
for
\begin{eqnarray*}
 &{}& \ba=\left(\tau_{11}\phi_1 - \tau_{12}\phi_2\right) \vect{v}_1
  +\left(\tau_{12}\phi_1 - \tau_{22}\phi_2\right) \vect{v}_2, \\
  &{}& \bb=\left(\tau_{11}\phi_1 - \tau_{12}\phi_2\right) \vect{v}_1
  +\left(\tau_{12}\phi_1 - \tau_{22}\phi_2\right) \vect{v}_2.
  \end{eqnarray*}
Substituting the following equations
\[
  \inner{\vect{v}_1}{\vect{v}_1}=\tau_{22},\qquad
  \inner{\vect{v}_1}{\vect{v}_2}=-\tau_{12},\qquad
  \inner{\vect{v}_2}{\vect{v}_2}=\tau_{11}
\]
into \eqref{eq:proof-MT}, 
we have that \eqref{eq:proof-MT} is equivalent to
\[
  -\Delta(\tau_{11}\phi_1^2 -2 \tau_{12} \phi_1\phi_2 + \tau_{22}\phi_2^2)=0.
\]
Since $\Delta>0$, we obtain \eqref{eq:MTS-eq}.
\enD

As a corollary of Theorem \ref{thm:graph},
we have the following.

\begin{Pro}\label{cor:L-flat-SMT}
A spacelike surface $M$ is strongly marginally trapped and 
totally $(\bn^T,\bn^S)$-umbilical if and only if
up to rigid motions of $\R^4_1,$ $M$ is given by
\begin{equation}\label{eq:Graph-f}
  \bX_f(u_1,u_2)=(f(u_1,u_2),f(u_1,u_2),u_1,u_2)),
\end{equation}
for a harmonic function $f:\R^2\rightarrow\R.$ 
\end{Pro}

\demo
By Proposition \ref{cor:L-flat-MT}, 
a totally $(\bn^T,\bn^S)$-umbilical marginally trapped surface
is given by $\bX_f$ as in \eqref{eq:Graph-f}
for a smooth function $f:\R^2\rightarrow\R$. 
Then,
the functions $\phi_1,~\phi_2$ 
defined by \eqref{eq:alpha} and \eqref{eq:beta}
can be calculated as
\[
  \phi_1=\phi_2=f_{u_1u_1}+f_{u_2u_2}.
\]
Thus, we have that $\bX_f$ 
is strongly marginally trapped if and only if
$f$ is harmonic.
\enD

\begin{Rem}
{\rm
The Classical Berntein theorem \cite{Be} says that
any minimal graph in $\mathbb{R}^3$ defined on the whole plane must be a plane.
By Proposition \ref{cor:L-flat-SMT},
we have that Bernstein-type theorem for strongly marginally trapped surfaces 
does not holds, that is, 
a strongly marginally trapped graph in $\mathbb{R}^4_1$
define on the whole plane does not need to be a spacelike plane.
}
\end{Rem}

%%%%%%%%%%%%%%%%%%%%%%%%
\section{Special cases}
In this section we naturally interpret minimal surfaces in Euclidean $3$-space, maximal surfaces in Lorentz-Minkowski $3$-space,
CMC$\pm 1$ surfaces in Hyperbolic $3$-space and CMC$\pm 1$ spacelike surfaces in de Sitter $3$-space as marginally trapped surfaces.
These surfaces have been well-investigated in the previous works.
We also show that intrinsic flat spacelike surfaces in the 
lightcone can be interpreted as marginally trapped surfaces.
\subsection{Surfaces in Euclidean space}
Let  $\mbox{\boldmath$X$}: U\longrightarrow \mathbb{R}^3$ be a surface in the Euclidean space
${\mathbb R}^3_0=\{\mbox{\boldmath$x$}\in \mathbb{R}^4_1\ |\ x_0=0\ \}.$ Then we can take
$\mbox{\boldmath$n$}^T=\mbox{\boldmath$e$}_0=(1,0,0,0)$, and have
\[
\mbox{\boldmath$n$}^S(u)= \frac{\mbox{\boldmath$e$}_0 \wedge \mbox{\boldmath$X$}_{u_1}(u) \wedge 
\mbox{\boldmath$X$}_{u_2}(u)}{\Vert \mbox{\boldmath$e$}_0 \wedge \mbox{\boldmath$X$}_{u_1}(u) \wedge 
\mbox{\boldmath$X$}_{u_2}(u) \Vert _1}\in S^2 \subset \mathbb{R}^3_0.
\]
Therefore, $\mbox{\boldmath$n$}^S(u)$ is the Euclidean unit normal of $M=\mbox{\boldmath$X$}(U)\subset\mathbb{R}^3_0$ at $p=\mbox{\boldmath$X$}(u).$
In this case, the lightcone Gauss map is given by
$\mathbb{L}^\pm (u)=\mbox{\boldmath$e$}_0\pm \mbox{\boldmath$n$}^S(u).$
So, the lightcone shape operator
is $\widetilde{S}^\pm _p=S_p(\be_0,\pm\bn^S)=-d(\mbox{\boldmath$e$}_0\pm \mbox{\boldmath$n$}^S)(u)
=\mp d\mbox{\boldmath$n$}^S(u)$ which is the Weingarten map for the surface in
the Euclidean space.
It follows that
$\widetilde{K}^\pm_\ell (u)=K_\ell (\be_0,\pm\bn^S(u))(u)=K(u)$ (i.e., the Gauss curvature) and
$\widetilde{H}^\pm _\ell (u)=H_\ell (\be_0,\pm\bn^S(u))(u)=\pm H(u)$ (i.e., the mean curvature).
Therefore, a lightcone flat spacelike surface is a developable surface and a 
marginally trapped surface is a minimal surface in this case.
Therefore, a marginally trapped surface is always a strongly marginally trapped surface.
\par
We remark that if $\mbox{\boldmath$n$}^T(u)=\mbox{\boldmath$v$}$ is a constant
timelike unit vector, the spacelike surface $M$ is a surface in the spacelike hyperplane $HP(\mbox{\boldmath$v$},c).$
Since $HP(\mbox{\boldmath$v$},c)$ is isometric to the Euclidean space $\mathbb{R}^3_0$, all the results for
$\mathbb{R}^3_0$ hold in this case.
\subsection{Spacelike surfaces in Minkowski $3$-space}
Let  $\mbox{\boldmath$X$}: U\longrightarrow \mathbb{R}^3_1$ be a spacelike surface in the Minkowski space
${\mathbb R}^3_1=\{\mbox{\boldmath$x$}=(x_0,x_1,x_2,x_3)\in \mathbb{R}^4_1\ |\ x_3=0\ \}.$ Then we can choose
$\mbox{\boldmath$n$}^S=\mbox{\boldmath$e$}_3=(0,0,0,1)$, and have
\[
\mbox{\boldmath$n$}^T(u)= \frac{\mbox{\boldmath$X$}_{u_1}(u) \wedge 
\mbox{\boldmath$X$}_{u_2}(u)\wedge \mbox{\boldmath$e$}_3}{\Vert \mbox{\boldmath$X$}_{u_1}(u) \wedge 
\mbox{\boldmath$X$}_{u_2}(u)\wedge \mbox{\boldmath$e$}_3 \Vert _1}\in H^2(-1) \subset \mathbb{R}^3_1.
\]
Therefore, $\mbox{\boldmath$n$}^T(u)$ is the timelike unit normal of $M=\mbox{\boldmath$X$}(U)\subset\mathbb{R}^3_1$ at $p=\mbox{\boldmath$X$}(u).$
In this case, the lightcone Gauss map is given by
$\mathbb{L}^\pm (u)=\mbox{\boldmath$n$}^T(u)\pm \mbox{\boldmath$e$}_3.$
So, the lightcone shape operator
is $\widetilde{S}^\pm _p=S_p(\bn^T(u),\pm\be _3)(u)=-d(\mbox{\boldmath$n$}^T(u)\pm \mbox{\boldmath$e$}_3)
=-d\mbox{\boldmath$n$}^T(u)$ which is the {\it spacelike shape operator} for the spacelike surface in
the Minkowski space.
It follows that
$\widetilde{K}^\pm_\ell (u)=K(\bn^T(u),\pm\be_3)(u)=K(u)$ (i.e., the Gauss-Kronecker curvature) and
$\widetilde{H}^\pm _\ell (u)=H_\ell (\bn^T(u),\pm\be _3)(u)=H(u)$ (i.e., the mean curvature).
Therefore, a lightcone flat spacelike surface is a spacelike developable surface and a
marginally trapped surface is a maximal surface.
In this case a marginally trapped surface is a strongly marginally trapped surface.
Of course, if $\bn^S(u)=\bv$ is constant spacelike unit vector, the spacelike surface $M$ is  surface
in the timelike hyperplane $HP(\mbox{\boldmath$v$},c).$
Since $HP(\mbox{\boldmath$v$},c)$ is isometric to the Minkowski space $\mathbb{R}^3_1$, all the results for
$\mathbb{R}^3_1$ hold in this case.
\subsection{Surfaces in the hyperbolic $3$-space}
Let  $\mbox{\boldmath$X$}: U\longrightarrow H^3_+(-1)$ be a surface in the
hyperbolic space. Then we can take
$\mbox{\boldmath$n$}^T=\mbox{\boldmath$X$}$, and
$$
\mbox{\boldmath$n$}^S(u)= \frac{\mbox{\boldmath$X$}(u) \wedge \mbox{\boldmath$X$}_{u_1}(u) \wedge 
\mbox{\boldmath$X$}_{u_2}(u)}{\Vert \mbox{\boldmath$X$}(u) \wedge \mbox{\boldmath$X$}_{u_1}(u) \wedge 
\mbox{\boldmath$X$}_{u_2}(u) \Vert _1}\in S^3_1 
$$
is univocally defined. We denote $\mbox{\boldmath$n$}^S$ by $\mbox{\boldmath$e$}$ and  call it the \textit{de Sitter normal vector field along $M$}.
The \textit{de Sitter Gauss image} is defined as the map
${\mathbb E}: U \rightarrow S^3_1$ given by ${\mathbb E}(u)= \mbox{\boldmath$e$} (u).$

We also have the \textit{hyperbolic Gauss image}
\[
\begin{array}{ccccl}
{\mathbb L}^{\pm} & : & U  & \longrightarrow & LC^*_+ \vspace{3mm}\\
& & u & \longmapsto & \mbox{\boldmath$X$}(u) \pm \mathbb{E}(u).
\end{array}
\]
In this particular case of surfaces in $\mathbb{R}^4_1$, the lightcone Gauss map on $M$ coincides with the {\it
hyperbolic Gauss map}, introduced in \cite{izu03} which is given by
$$
\widetilde{{\mathbb L}}^\pm(u)=\widetilde{ \mbox{\boldmath$X$}(u) \pm {\mathbb E}(u)}.
$$ 

We observe that this notion of hyperbolic Gauss map is equivalent to the one defined \cite{Bry,Ep,Koba1,Koba2} in other models.
\par
It is well-known that the hyperbolic space is a model for
the non-Euclidean geometry of Gauss-Bolyai-Lobachevsky (i.e., the hyperbolic geometry).
Recently, it has been discovered another geometry in the hyperbolic space which is called the horospherical geometry
(cf., \cite{izu03,iz05,iz06,IST10}).
In this case
$K_\ell (\bX(u),\mathbb{E}(u))=K^\pm _h(u)$ (i.e., the {\it hyperbolic Gauss-Kronecker curvature}
\cite{izu03}),
$\widetilde{K}^\pm_\ell (u)=\widetilde{K}^\pm_h(u)$ (i.e., the {\it horospherical Gauss-Kronecker curvature} \cite{iz06}) and
$H_\ell (\bX(u),\pm\mathbb{E}(u)) (u)=H^\pm _h(u)$, where $H^\pm _h=-1\pm H$ is the {\it hyperbolic mean curvature} \cite{izu03} and $H$ is the mean curvature in the ordinary sense (i.e., the {\it de Sitter mean curvature} \cite{izu03}).
Therefore, a lightcone flat spacelike surface is a horospherical flat surface\cite{IST10} and a
marginally trapped surface is a CMC$\pm 1$ surface.
Since $\bn^T=\bX,$ $H(\bn^T)\equiv -1$, there are no strongly marginally trapped surfaces in the 
hyperbolic space.
Let $V=\alpha\mathbb{L}^\pm$ be the variation vector filed.
Then the first variational formula is
\[
 \left. \frac{d}{d \ve} \right|_{\ve=0} A(\bX^{\ve})=-2\int_{U} \alpha H^\pm_h dM
=-2\int_{U} \alpha (-1\pm H)dM
\]
and the second variational formula is
\[
  \left. \frac{d^2}{d \varepsilon^2} \right|_{\varepsilon=0} A(\bX^{\varepsilon})
  = 2 \int_{U} \alpha ^2K^\pm_h dM.
\]

\subsection{Spacelike surfaces in de Sitter $3$-space}

Let  $\mbox{\boldmath$X$}: U\longrightarrow S^3_1$ be a spacelike surface in the
de Sitter $3$-space. Then we can take
$\mbox{\boldmath$n$}^S=\mbox{\boldmath$X$}$, and
$$
\mbox{\boldmath$n$}^T(u)= \frac{\mbox{\boldmath$X$}_{u_1}(u) \wedge 
\mbox{\boldmath$X$}_{u_2}(u)\wedge \mbox{\boldmath$X$}(u) }{\Vert \mbox{\boldmath$X$}_{u_1}(u) \wedge 
\mbox{\boldmath$X$}_{u_2}(u)\wedge \mbox{\boldmath$X$}(u)  \Vert _1}\in H^3(-1) 
$$
is defined. We denote $\mbox{\boldmath$n$}^S$ by $\mbox{\boldmath$e$}$ and  call it the \textit{hyperbolic normal vector field along $M$}.
The \textit{hyperbolic Gauss image} is defined as the map
${\mathbb E}: U \rightarrow H^3(-1)$ given by ${\mathbb E}(u)= \mbox{\boldmath$e$} (u).$

We also have the \textit{lightcone Gauss image}
\[
\begin{array}{ccccl}
{\mathbb L}^{\pm} & : & U  & \longrightarrow & LC^*_+ \vspace{3mm}\\
& & u & \longmapsto & \mbox{\boldmath$X$}(u) \pm \mathbb{E}(u).
\end{array}
\]
In this particular case of spacelike surfaces in $\mathbb{R}^4_1$, the lightcone Gauss image on $M$ was introduced in \cite{K09}.
We only give the following proposition as a corollary of \cite[Proposition 4.5]{IR07}.
\begin{Pro}
Let  $\mbox{\boldmath$X$}: U\longrightarrow S^3_1$ be a spacelike surface in the
de Sitter $3$-space. Then the lightcone Gauss image ${\mathbb L}^\pm(u) $ is
a constant vector \mbox{\boldmath$v$} if and only if
$M$ is a parabola defined by $HP(\mbox{\boldmath$v$},c)\cap S^3_1.$
\end{Pro}
The parabola $HP(\mbox{\boldmath$v$},c)\cap S^3_1$ for a lightlike vector \mbox{\boldmath$v$} is
called a {\it de Sitter horosphere.}
The geometry related to the hyperbolic Gauss image might be called {\it de Sitter horospherical geometry}.
Further results on de Sitter horospherical geometry are referred in the articles \cite{K09,K10}.
\par
A marginally trapped surface in the de Sitter space is also a CMC$\pm 1$ surface and there are
no strongly marginally trapped surfaces.
In this case, we have the similar formulae for the variations to the case of surfaces in Hyperbolic space.

\subsection{Spacelike surfaces in the lightcone}

The induced metric is degenerate on the lightcone, so that the ordinary submanifolds theory cannot work for
surfaces in the lightcone.
In \cite{I09} we constructed the basic tools for the study of the extrinsic geometry on spacelike surfaces 
in the lightcone $LC^*$ (see, also \cite{L07}) as one of the applications of the mandala of Legendrian dualitites 
between pseudo-spheres in Minkowski space-time\cite{I09,IY11}.
We define one-forms $\langle d\bv,\bw\rangle _1 =-w_0dv_0+\sum_{i=1}^3 w_idv_i$,
$\langle \bv,d\bw \rangle _1 =-v_0dw_0+\sum_{i=1}^3 v_idw_i$ on $\R^{4}_1\times\R^{4}_1$ and consider the following double fibration with one-forms
\par
(a) $LC^*\times LC^*\supset \Delta _4=\{(\bv ,\bw )\ |\ \langle \bv ,\bw \rangle _1 =-2\ \}$,
\par
(b) $\pi _{41}:\Delta _4\lon LC^*$,$\pi _{42}:\Delta _4\lon LC^*$,
\par
(c) $\theta _{41}=\langle d\bv ,\bw\rangle _1 |\Delta _4$,
$\theta _{42}=\langle \bv ,d\bw\rangle _1 |\Delta _4$,
\par\noindent
where $\pi _{41}(\bv,\bw)=\bv$, $\pi _{42}(\bv ,\bw)=\bw$ are the canonical projections.
Moreover, $\theta _{41}=\langle d\bv ,\bw\rangle _1|\Delta _4$ and $\theta _{42}=\langle d\bw,\bv\rangle _1 |\Delta _4$
are the restrictions of the one-forms $\langle d\bv ,\bw\rangle _1$ and $\langle d\bw ,\bv\rangle _1$
on $\Delta _4.$ It has been shown in \cite{I09} that $(\Delta _4,K_4)$ is a contact manifold and
both of $\pi _{4j}$ $(j=1,2)$ are Legendrian fibrations, where $K_4=\theta _{41}^{-1}(0)=\theta _{42}^{-1}(0).$
In \cite{I09} we defined four Legendrian fibrations $(\Delta _i,K_i)$ $(i=1,2,3,4)$ such that these are contact diffeomorphic 
to each other. Here, we only use $(\Delta _4,K_4).$
For definitions and basic results of Legendrian fibrations, see \cite{Arnold1,I09}.
\par
Let $\bX:U\lon LC^*$ be a spacelike surface in $LC^*\subset \R^4_1.$
In \cite{I09} it has been shown that there exists a unique map $\bX^\ell :U\lon LC^*$ such that
$\mathscr{L}_4:U\lon \Delta _4$, defined by $\mathscr{L}_4(u)=(\bX(u), \bX^\ell (u))$, is a Legendrian embedding.
We call $\bX^\ell$ the {\it lightcone Gauss image} of $M=\bX(U).$
Applying the basic properties of $\mathscr{L}_4$ as a Legendrian embedding , we defined curvatures 
of $M=\bX(U)$ in \cite{I09} by using $\bX^\ell$ as a normal vector field.
The lightcone shape operator is defined to be $S^\ell _p=-d\bX^\ell (u)$ for
$p=\bX(u).$
The {\it lightcone Gauss-Kronecker curvature} is $K_\ell (p)=\det S^\ell_p$ and
the {\it lightcone means curvature} is $H_\ell (p)=({\rm Trace}\,S^\ell_p)/2.$
Since $\langle \bX(u),\bX^\ell (u)\rangle _1 =-2,$  we have
\[
\bn^T(u)=\frac{\bX(u)+\bX^\ell (u)}{2}\ {\rm and}\
\bn^S(u)=\frac{\bX(u)-\bX^\ell (u)}{2}.
\]
It follows that
\[
\kappa _i(\bn^T)(p)=\frac{1}{2}(-1+\kappa ^\ell _i(p)), \ \kappa _i(\bn^S)(p)=\frac{1}{2}(-1-\kappa ^\ell _i(p)),
\]
where $\kappa ^\ell _i(p)$ $(i=1,2)$ are the lightcone principal curvature of $M=\bX(U)$ at $p=\bX(u)$ (i.e., the eigenvalues of 
$S^\ell _p$). 
One of the interesting properties of these curvatures is the following \lq\lq Theorema Egregium\rq\rq (\cite{I09,L07}).
\begin{Th}[Theorem 9.3 in \cite{I09}]
Let $K_I$ be the intrinsic Gauss curvature of $M=\bX(U).$ Then we have the relation
\[
K_I=H_\ell.
\]
\end{Th}
\par
It follows that a marginally trapped surface is an intrinsic flat surface as the
following theorem shows.
\begin{Th}
Let $\bX:U\lon LC^*$ be a spacelike surface and $\bX^\ell$ the $\Delta _4$-dual of $\bX$.
Then the following conditions are equivalent{\rm:}
\par\noindent
{\rm 1)} $M=\bX(U)$ is a marginally trapped surface.
\par\noindent
{\rm 2)} $H_\ell\equiv 0.$
\par\noindent
{\rm 3)} The mean curvature vector $\mathfrak{H}(p)$ is zero or parallel to $\bX(u)$
at any $p=\bX(u).$
\par\noindent
{\rm 4)} $K_I\equiv 0.$
\par\noindent
{\rm 5)} The Gauss curvature vector $\mathfrak{K}(p)$ is zero or parallel to  $\bX(u)$
at any $p=\bX(u).$
\end{Th}
\demo
Since 
\[
\kappa _i(\bn^T)(p)=\frac{1}{2}(-1+\kappa ^\ell _i(p)), \ \kappa _i(\bn^S)(p)=\frac{1}{2}(-1-\kappa ^\ell _i(p)),
\]
we have
\begin{eqnarray*}
&{}&H(\bn^T)(u)=\frac{1}{2}(-2+\kappa ^\ell _1(p)+\kappa ^\ell _2(p))=H_\ell(u) -1,\\
&{}&H(\bn^S)(u)=\frac{1}{2}(-2-\kappa ^\ell _1(p)-\kappa ^\ell _2(p))=-H_\ell(u) -1.
\end{eqnarray*}
It follows that
\begin{eqnarray*}
H_\ell (\bn^T,\bn^S)(u)&=&H(\bn^T)(u)+H(\bn^S)(u)=-2,\\
 H_\ell (\bn^T,-\bn^S)(u)&=& H(\bn^T)(u)-H(\bn^S)(u)=2H_\ell(u).
\end{eqnarray*}
Therefore, $H_\ell (u)=0$ if and only if $H(\bn^T)(u)=H(\bn^S)(u).$
This means that $\mathfrak{H}(p)$ is parallel to $\bX(u)=\bn^T(u)+\bn^S(u).$
By definition, the condition 3) implies the condition 1).
By Corollary 4.4, $M=\bX(U)$ is marginally trapped if and only if
$H_\ell \equiv 0.$
Since $H_\ell =K_I,$ the conditions 2) and 4) are equivalent.
By Proposition 3.13,
the condition 4) (i.e., $K_I(p)=0$) is equivalent to the condition that
$\mathfrak{K}(p)=\bo$ or $\mathfrak{K}(p)$ is parallel to $\bn^T(u)+\bn^S(u)=\bX(u).$
\enD
\par
Since $H_\ell (\bn^T,\bn^S)(u)=-2,$ there are no strongly marginally trapped surfaces in
the lightcone.
\par
Let $V=\alpha\bX^\ell$ be the variation vector filed.
Then the first variation formula is
\[
\left. \frac{d}{d \ve} \right|_{\ve=0} A(\bX^{\ve})=-2\int_{U} \alpha H_\ell dM
=-2 \int_{U} \alpha K_I dM
\]
and the second variation formula is
\[
  \left. \frac{d^2}{d \varepsilon^2} \right|_{\varepsilon=0} A(\bX^{\varepsilon})
  = 2 \int_{U} \alpha ^2K_\ell  dM.
\]


\begin{thebibliography}{1}

\bibitem{AGM05}
Aledo, J. A., G\'alvez, J. A. and Mira, P.
{\it Marginally Trapped Surfaces in $\mathbb{L}^4$ and an Extended Weierstrass-Bryant Representaion}.
Annals of Global Analysis and Geometry {\bf 28}, 395--415 (2005)

\bibitem{AM}
Andersson, L. and Metzger, J.
{\it Curvature estimates for stable marginally trapped surfaces},
J. Differential Geom. {\bf 84}, 231--265 (2010)

\bibitem{Arnold1}
Arnol'd, V. I., Gusein-Zade, S. M. and Varchenko, A. N.
{\it Singularities of Differentiable Maps vol. I}.
Birkh\"auser (1986) 

\bibitem{Be}
Bernstein, S.
{\it Sur un th\'eor\`eme de g\'eom\'etrie et ses applications aux \'equations aux d\'eriv\'ees partielles du type elliptique}.  
Comm. de la Soc. Math. de Kharkov (2\'eme s\'er.) {\bf 15}, 38--45 (1915--1917)
  
\bibitem{Bry}
Bryant, R. L.
{\it Surfaces of mean curvature one in hyperbolic space}.
in Th\'eorie des vari\'et\'es minimales et applications (Palaiseau, 1983--1984),
Ast\'erisque No. 154--155 (1987), 12, 321--347, 353 (1988) 
\bibitem{ChI91}
Chen, B. Y. and Ishikawa, S.
{\it Biharmonic surfaces in pseudo-Euclidean spaces}.
Memoirs Fac. Sci. Kyushu Univ. Ser. A, Math. {\bf 45}, 325--349 (1991)
\bibitem{Ch09} Chen, B. Y.
{\it Marginally trapped surfaces and Kaluza-Klein theory}.
International Electric Journal of Geometry {\bf 2}, 1--16 (2009)

\bibitem{Cher}
Chern, S. S.
{\it An elementary proof of the existence of isothermal parameters on a surface}.
 Proceedings of AMS {\bf 6}, 771--782 (1955)

\bibitem{CGP10}
Chru\'sciel, P. T., Galloway, G. J. and Pollack, D.
{\it Mathematical general relativity: a sampler}.
Bulletin of the AMS. {\bf 47}, 567--638 (2010)

\bibitem{Ep} Epstein, C. L.
{\it The hyperbolic Gauss map and quasiconformal reflections}.
J. Reine Angew. Math. {\bf 372}, 96--135 (1986)

\bibitem{izu03}
Izumiya, S., Pei, D-H. and Sano, T. 
{\it Singularities of hyperbolic
Gauss maps}. 
Proc. London Math. Soc.  {\bf 86}, 485--512 (2003)

\bibitem{iz05} Izumiya, S.  Pei, D-H. Romero-Fuster, M. C. and 
Takahashi, M. {\it Horospherical geometry of submanifolds in hyperbolic $n$-space}. 
Journal of London Mathematical
Society, \textbf{71}, 779--800 (2005)

\bibitem{iz06} Izumiya, S. and Romero Fuster, M. C.
{\it The horospherical Gauss-Bonnet type theorem in hyperbolic
  space}. Journal of Math. Soc. Japan \textbf{58},
965--984 (2006)
\bibitem{IKPR06}
Izumiya, S., Kossowski, M., Pei, D-H. and Romero Fuster, M. C.
{\it Singularities of lightlike hypersurfaces in Minkowski
$4$-space}. Tohoku Math. J. (2) {\bf 58}, 71--88 (2006)

\bibitem{IR07}
Izumiya, S. and Romero Fuster, M. C.
{\it The lightlike flat geometry on spacelike submanifolds of codimension two
in Minkowski space}. Selecta Mathematica (NS) {\bf 13} 23--55 (2007)
\bibitem{IRS09}
Izumiya, S., Romero Fuster, M. C. and Saji, K.
{\it Flat lightlike hypersurfaces in Lorentz-Minkowski 4-space}.
J. of Geometry and Physics {\bf 59}, 1528--1546 (2009)
\bibitem{I09}
Izumiya, S.
{\it Legendrian dualities and spacelike hypersurfaces in the lightcone}.
 Moscow Mathematical Journal, {\bf 9}, 325--357 (2009)
\bibitem{IST10}
Izumiya, S., Saji, K and Takahashi, M.
  {\it Horospherical flat surfaces in Hyperbolic
        $3$-space}.
  J. Math. Soc. Japan {\bf 62}, No.3, 789-849 (2010)
  
\bibitem{IY11} Izumiya, S. and Y\i ld\i r\i m, H.
{\it
Extensions of the mandala of Legendrian dualities for pseudo-spheres in Lorentz-Minkowski space}.
Topology and its Applications  {\bf 159}, 509--518 (2012)
  
\bibitem{K09} Kasedou, M. {\it Singularities of lightcone Gauss
images of spacelike hypersurfaces in de Sitter space}. Journal of Geometry, \textbf{94},  107--121
(2009) 
\bibitem{K10} Kasedou, M. {\it Spacelike submanifolds of codimension two in de Sitter space}. Journal of Geometry and Physic, \textbf{60}, 31--42 (2010)
 
\bibitem{Koba1}
Kobayashi, T.
{\it Null varieties for convex domains {\rm (}Japanese{\rm )}}.
Reports on unitary representation seminar {\bf 6}, 1--18 (1986)
\bibitem{Koba2}
Kobayashi, T.
{\it Asymptotic behaviour of the null variety for a convex domain in a
non-positively curved space form}.
Journal of the Faculty of Science, University of Tokyo, {\bf 36}, 389--478 (1986)


\bibitem{L07} Liu, H. {\it Surfaces in the lightlike cone}. Journal
of Math. Annal. Appl. \textbf{325}, 1171--1181 (2007)


\bibitem{Oneil}
O'Neill, B.
{\it Semi-Riemannian Geometry},
Academic Press, New York (1983)
\bibitem{Pen65}
Penrose, R.
{\it Gravitaional collapse and space-time singularities}.
Phys. Rev. Lett. {\bf 14}, 57--59 (1965)

\end{thebibliography}
\end{document}